\documentclass[12pt]{amsart}
\usepackage{amsbsy}
\usepackage{graphicx,epsfig,subfigure,psfrag}
\begin{document}

\newtheorem{theorem}{Theorem}
\newtheorem{proposition}{Proposition}
\newtheorem{lemma}{Lemma}
\newtheorem{corollary}{Corollary}
\newtheorem{definition}{Definition}
\newtheorem{remark}{Remark}
\newcommand{\tex}{\textstyle}
\numberwithin{equation}{section} \numberwithin{theorem}{section}
\numberwithin{proposition}{section} \numberwithin{lemma}{section}
\numberwithin{corollary}{section}
\numberwithin{definition}{section} \numberwithin{remark}{section}
\newcommand{\ren}{\mathbb{R}^N}
\newcommand{\re}{\mathbb{R}}
\newcommand{\n}{\nabla}
\newcommand{\iy}{\infty}
\newcommand{\pa}{\partial}
\newcommand{\fp}{\noindent}
\newcommand{\ms}{\medskip\vskip-.1cm}
\newcommand{\mpb}{\medskip}
\newcommand{\AAA}{{\bf A}}
\newcommand{\BB}{{\bf B}}
\newcommand{\CC}{{\bf C}}
\newcommand{\DD}{{\bf D}}
\newcommand{\EE}{{\bf E}}
\newcommand{\FF}{{\bf F}}
\newcommand{\GG}{{\bf G}}
\newcommand{\oo}{{\mathbf \omega}}
\newcommand{\Am}{{\bf A}_{2m}}
\newcommand{\CCC}{{\mathbf  C}}
\newcommand{\II}{{\mathrm{Im}}\,}
\newcommand{\RR}{{\mathrm{Re}}\,}
\newcommand{\eee}{{\mathrm  e}}
\newcommand{\LL}{L^2_\rho(\ren)}
\newcommand{\LLL}{L^2_{\rho^*}(\ren)}
\renewcommand{\a}{\alpha}
\renewcommand{\b}{\beta}
\newcommand{\g}{\gamma}
\newcommand{\G}{\Gamma}
\renewcommand{\d}{\delta}
\newcommand{\D}{\Delta}
\newcommand{\e}{\varepsilon}
\newcommand{\var}{\varphi}
\renewcommand{\l}{\lambda}
\renewcommand{\o}{\omega}
\renewcommand{\O}{\Omega}
\newcommand{\s}{\sigma}
\renewcommand{\t}{\tau}
\renewcommand{\th}{\theta}
\newcommand{\z}{\zeta}
\newcommand{\wx}{\widetilde x}
\newcommand{\wt}{\widetilde t}
\newcommand{\noi}{\noindent}
\newcommand{\uu}{{\bf u}}
\newcommand{\xx}{{\bf x}}
\newcommand{\yy}{{\bf y}}
\newcommand{\zz}{{\bf z}}
\newcommand{\aaa}{{\bf a}}
\newcommand{\cc}{{\bf c}}
\newcommand{\jj}{{\bf j}}
\newcommand{\ggg}{{\bf g}}
\newcommand{\UU}{{\bf U}}
\newcommand{\YY}{{\bf Y}}
\newcommand{\HH}{{\bf H}}
\newcommand{\GGG}{{\bf G}}
\newcommand{\VV}{{\bf V}}
\newcommand{\ww}{{\bf w}}
\newcommand{\vv}{{\bf v}}
\newcommand{\hh}{{\bf h}}
\newcommand{\di}{{\rm div}\,}
\newcommand{\ii}{{\rm i}\,}
\newcommand{\inA}{\quad \mbox{in} \quad \ren \times \re_+}
\newcommand{\inB}{\quad \mbox{in} \quad}
\newcommand{\inC}{\quad \mbox{in} \quad \re \times \re_+}
\newcommand{\inD}{\quad \mbox{in} \quad \re}
\newcommand{\forA}{\quad \mbox{for} \quad}
\newcommand{\whereA}{,\quad \mbox{where} \quad}
\newcommand{\asA}{\quad \mbox{as} \quad}
\newcommand{\andA}{\quad \mbox{and} \quad}
\newcommand{\withA}{,\quad \mbox{with} \quad}
\newcommand{\orA}{,\quad \mbox{or} \quad}
\newcommand{\ef}{\eqref}
\newcommand{\ssk}{\smallskip}
\newcommand{\LongA}{\quad \Longrightarrow \quad}
\def\com#1{\fbox{\parbox{6in}{\texttt{#1}}}}
\def\N{{\mathbb N}}
\def\A{{\cal A}}
\newcommand{\de}{\,d}
\newcommand{\eps}{\varepsilon}
\newcommand{\be}{\begin{equation}}
\newcommand{\ee}{\end{equation}}
\newcommand{\spt}{{\mbox spt}}
\newcommand{\ind}{{\mbox ind}}
\newcommand{\supp}{{\mbox supp}}
\newcommand{\dip}{\displaystyle}
\newcommand{\prt}{\partial}
\renewcommand{\theequation}{\thesection.\arabic{equation}}
\renewcommand{\baselinestretch}{1.1}
\newcommand{\Dm}{(-\D)^m}

\title
 {\bf Incomplete  self-similar blow-up in a semilinear
  fourth-order
 reaction-diffusion
 equation}

\author {V.A.~Galaktionov}

\address{Department of Mathematical Sciences, University of Bath,
 Bath BA2 7AY, UK}
\email{vag@maths.bath.ac.uk}



\keywords{4th-order semilinear parabolic equation, incomplete
blow-up, self-similar solutions, similarity extensions.}

 \subjclass{35K55, 35K40  }
\date{\today}

\begin{abstract}

Blow-up behaviour
 for the
4th-order semilinear reaction-diffusion  equation
  \be
  \label{01}
u_t= - u_{xxxx} + |u|^{p-1} u \inB \re \times \re_+, \,\, p>1,
  \ee
is studied. For the classic semilinear heat equation from
combustion theory
 \be
 \label{02}
 u_t= u_{xx}+ u^p \quad (u \ge 0),
  \ee
various blow-up patterns were investigated
 since 1970s, while the case of higher-order diffusion was
studied much less. Blow-up self-similar solutions of \ef{01},
 $$
  \tex{
 u_{-}(x,t)=(T-t)^{-\frac 1{p-1}} f(y), \quad y= \frac
 x{(T-t)^{1/4}},
 }
 $$
are shown to admit global extensions for $t>T$ in an analogous
similarity form:
 $$
  \tex{
 u_{+}(x,t)=(t-T)^{-\frac 1{p-1}} F(y), \quad y= \frac
 x{(t-T)^{1/4}}.
 }
 $$
  The continuity at $t=T$ is preserved in the sense that
 $$
  u_{-}(x,T^-)=u_+(x,T^+)=C_1|x|^{-\frac 4{p-1}} \,\,\, \mbox{for all} \,\,\, x \in \re\setminus\{0\}
   \quad (C_1={\rm const.} \ne
  0).
   $$
 This is in a striking
difference with blow-up for \ef{02}, which is known to be always
{\em complete} in the sense that the minimal (proper) extension
beyond blow-up is  $u(x,t) \equiv + \iy$ for $t>T$. Difficult
4th-order dynamical systems for $f$ and $F$ are studied by a
combination of various analytic, formal, and careful numerical
methods.
 Other non-similarity  patterns
for \ef{01} with non-generic   {\em complete blow-up} are also
discussed.

\end{abstract}

\maketitle

\section{Introduction: self-similar blow-up patterns of higher-order reaction-diffusion equations
and main results}
  \label{S.1}

 \subsection{Fourth-order RDE and blow-up}

This paper continues the study began in \cite{BGW1, Gal2m,
GalBlow5, GW1} of blow-up patterns
 for the
{\em fourth-order reaction-diffusion equation} (the RDE--4)
\be
\label{m2}
  \tex{
u_t= - u_{xxxx} + |u|^{p-1} u \inB \re \times \re_+ \whereA p>1.
 }
 \ee
 For
 applications of such higher-diffusion models, see  surveys and references in
 \cite{BGW1, GW1}. In general,
 higher-order semilinear parabolic equations arise in many physical
applications such as thin film theory, convection-explosion
theory, lubrication theory, flame and wave propagation (the
Kuramoto-Sivashinsky equation and the extended Fisher-Kolmogorov
equation), phase transition at critical Lifschitz points,
bi-stable systems and applications to structural mechanics; the
effect of fourth-order terms on self-focusing problems in
nonlinear optics are also well-known in applied and mathematical
literature. For a systematic treatment of extended KPPF-equations;
see in Peletier--Troy  \cite{PelTroy}.

 The RDE--4 \ef{m2} has
  the
 bi-harmonic diffusion operator  $-D_x^4$ and is a higher-order counterpart
 of classic semilinear heat equation from combustion theory \cite{ZBLM}
\be
  \label{m1}
u_t= u_{xx} + |u|^{p-1}u \inB \re \times \re_+,
 \ee
 with many well-known properties including various aspects of
 blow-up behaviour; see
 a number of well-known monographs
 \cite{BebEb, SGKM, GSVR, Pao, MitPoh, AMGV, GalGeom, QSupl}.
Surveys in \cite{BGW1} and in a more recent paper \cite{GalBlow5}
 contain necessary information concerning relations between these
 models \ef{m2} and \ef{m1} and description of similarity and
 other blow-up patterns for \ef{m2}.

Note  that another related fourth-order one-dimensional semilinear
parabolic equation
\begin{equation}
\label{JS}
 u_t =
   - u_{xxxx} - [(2-(u_x)^2)u_x]_x - \a u + qe^{su},
\end{equation}
where $\a $, $q$ and $s$ are positive constants obtained from
physical parameters,
 occurs in the Semenov-Rayleigh-Benard problem
\cite{JMS}, where the equation is derived  in studying the
interaction between natural convection and the explosion of an
exothermically-reacting fluid confined between two isothermal
horizontal plates. This is an evolution equation for the
temperature fluctuations in the presence  of natural convection,
wall losses and chemistry. It can be considered as a formal
combination of the equation derived in \cite{GS} (see also
\cite{CP}) for the Rayleigh-Benard problem and of the Semenov-like
energy balance \cite{Sem, F-K} showing that  natural convection
and the explosion mechanism may reinforce each other;
 see more details on physics
and mathematics of blow-up in \cite{GW1}.
 In a special limit, \ef{JS} reduces to the {\em generalized
 Frank-Kamenetskii equation}
 \begin{equation}
\label{u4e} u_t = -u_{xxxx} + e^u,
\end{equation}
 which plays a role of a natural extension of the classic
 Frank-Kamenetskii equation
  \be
   \label{FK11}
   u_t=\D u+{\mathrm e}^u,
    \ee
     derived  in
 solid fuel theory in the 1930s, \cite{Fr-K}.

Thus, similar to the second-order equations pair \ef{m1},
\ef{FK11}, where both have  an equal physical significance, we
choose the RDE--4 \ef{m2} as a leading model, though some of the
results are naturally applied to \ef{u4e}.
  Note that \ef{m2} can be considered as a
 non-mass-conservative counterpart of the well-known {\em limit unstable
 Cahn--Hilliard equation} from phase transition,
  \be
   \label{CH1}
   u_t= - u_{xxxx} - (|u|^{p-1}u)_{xx} \inB \re \times \re_+,
    \ee
 which is known to admit various families of blow-up solutions; see \cite{EGW1}
 for a long list of references.
 It is key that the mass-conservation for $L^1$-solutions of
 \ef{CH1} naturally demands existence of an extension after
 blow-up, which was studied in \cite{GalJMP}. This is not that
 straightforward for the current model \ef{m2}.
  Indeed, \ef{m2} is related to the
 famous mass-conserving divergence  {\em Kuramoto--Sivashinsky equation} from flame
 propagation theory
  \be
  \label{KS1}
  u_t= - u_{xxxx} -u_{xx} + u u_x \inB \re \times \re_+,
   \ee
 which always admits  global solutions, so no blow-up for \ef{KS1} exists.

 As usual in PDE theory, {\em blow-up} means that, in the
 Cauchy problem (the CP) for \ef{m2} (or \ef{m1}) with bounded smooth initial function $u_0(x)$,
the classic bounded solution $u=u(x,t)$   exists in $\re \times
(0,T)$, while
 \be
 \label{Bl1}
 \tex{
  \sup_{x \in \ren} \, |u(x,t)| \to +\iy \asA t \to T^-,
  }
  \ee
  where $T \in \re_+=(0,+\iy)$, depending on given data $u_0$, is
   called the {\em blow-up time} of
  the solution $u(x,t)$.
 It is also convenient to use the quite popular nowadays
 auxiliary
classification from Hamilton
 \cite{Ham95}, where Type I blow-up
 means the solutions satisfying, for some constant $C >0$ (depending
 on $u$),
  \be
  \label{S5}
  \begin{matrix}
 \mbox{Type I:} \quad
 (T-t)^{\frac 1{p-1}} |u(x,t)| \le C  \asA t \to T^-,
  \,\,\,\mbox{and, otherwise}\ssk\ssk\\
  \mbox{Type II:} \quad \limsup_{t \to T^-} (T-t)^{\frac 1{p-1}}
  \sup_x |u(x,t)|=+\iy \qquad\qquad\,\,\,
  \end{matrix}
   \ee
   (Type II also called {\em slow} blow-up in
   \cite{Ham95}).
   In R--D theory, blow-up with the first dimensional estimate in \ef{S5}
   was usually called of {\em self-similar} rate, while Type II was
   referred to as {\em fast} and {\rm non self-similar}; see
   \cite{AMGV} and \cite{SGKM}.

\subsection{Main results: on incomplete blow-up and extended global semigroup}

Unlike the previous papers, we are now more interested in a
possible extension of blow-up solutions beyond blow-up, i.e., for
$t>T$. For \ef{m1},  this is not possible: blow-up of any its
solution $u \ge 0$ is known to be {\em complete}, i.e., the proper
(minimal) extension is
 \be
 \label{ext1}
  u(x,t)=+\iy \inB \re \times (T,+\iy).
   \ee
 In fact, \ef{ext1} is guaranteed by the Maximum Principle; see
 \cite{GV} and \cite[Ch.~6,7]{GalGeom} for main concepts and results of extended semigroup theory and
 further references.

Our main goal is to justify, that, rather surprisingly for us, for
the RDE--4  (\ref{m2}),
 \be
 \label{ex2}
  \fbox{$
\mbox{self-similar blow-up is {\em incomplete} and admits
self-similar extension for $t>T$.}
 $}
 \ee
Actually, \ef{ex2} well corresponds to {\em Leray's scenario} of
self-similar blow-up and similarity extension beyond, proposed by
Jean Leray in 1934 for the Navier--Stokes equations in $\re^3$;
see \cite[p.~245]{Ler34} for the precise  formulation and
 \cite{GalJMP} for recent discussions.
For equation \ef{m2}, \ef{ex2} has an important corollary:
 \be
 \label{cor1}
 \fbox{$
  \mbox{the number of self-similar extensions for $t>T$ is not more than
  countable}
   $}
   \ee
   (actually, it is finite, but a proof exists in the analytic case $p=3,\,5\,...$ only).
   This confirms a plausible existence of an extended semigroup of unique
   global proper (``minimal") solutions of the Cauchy problem  for \ef{m2}.

We will also present some convincing facts that self-similar
blow-up is a generic (structurally stable) one for \ef{m2}, so
that \ef{ex2}  implies a certain  possibility of existence of a
kind of extended semigroup theory of  unique global proper
solutions of \ef{m2} defined for all $t>0$. For \ef{m1}, this is
well known; see \cite{GV} as a guide. For \ef{m2},
 this question will be also discussed but
essentially remains open.

  For \ef{m1}, the property \ef{ex2} is revealed for the RDE--2 in $\ren$
  only, i.e., for the equation
  \be
  \label{m1N}
   \tex{
  u_t= \D u + u^p \whereA p > p_{\rm Sob}= \frac{N+2}{N-2}, \quad
  N>2;
  }
  \ee
 see \cite{GV}. However, even for \ef{m1N}, such an incomplete blow-up
is not generic in the sense that, for almost all (a.a.)
 solutions,
blow-up is expected to be complete. This is not fully justified
rigorously, though the results imply this at least in the radial
geometry.

Thus, according to \ef{ex2},  the RDE--4 \ef{m2} exhibits more
flexibility in adapting blow-up solutions for $t>T$, though  the
analysis becomes essentially more difficult. To justify \ef{ex2}
in Section \ref{S2} (blow-up patterns for $t<T$) and \ref{S3}
(global patterns for $t>T$), we will use a variety of methods
including some analytical ones, but the final conclusions will
eventually depend on careful numerical experiments. This is an
unavoidable feature of the study of the 4th-order dynamical
systems to be derived, which seems do not admit a fully
mathematically rigorous investigation.
 As a key issue, we claim that
 \be
 \label{k1}
  \fbox{$
  \mbox{incomplete self-similar blow-up (\ref{ex2}) for (\ref{m2}) has a
  pure dimensional nature.}
 $}
  \ee
  Note that earlier, by a similar reason \ef{k1},
  existence of extended solutions was obtained \cite{GalJMP}
for the {
 Cahn--Hilliard equation} \ef{CH1} for $p=3$:
  \be
  \label{mm7}
  u_t=-u_{xxxx}-(u^3)_{xx} \inB \re \times \re_+,
   \ee
for which existence of a countable family of positive blow-up
similarity solutions was established  in \cite{EGW1}. However,
since \ef{mm7} is divergent and hence  conservative by preserving
the total mass of $L^1$-solutions, existence of  solutions beyond
blow-up is more natural than for the non-conservative and
non-divergent model \ef{m2}.

 In Section \ref{S4}, we present  other blow-up patterns for \ef{m2},
 introduced first in \cite{Gal2m}, which do not admit global
 extensions beyond blow-up, so belong to the case of {\em complete
 blow-up}.
 Fortunately, such complete blow-up is most plausibly non-generic for
 \ef{m2}.

We hope that the conclusion \ef{k1} and others will be helpful for
understanding singularity formation and extension concepts for
$2m$th-order nonlinear evolution PDEs including
 \be
 \label{mm98}
  u_t= - (-\D)^m u + |u|^{p-1}u \quad (m \ge 2),
   \ee
 which are
definitely short of new ideas concerning construction of extended
semigroups of unique global solutions. Note that some previous
results in \cite{BGW1, Gal2m} were already oriented to arbitrary
$m=2,3,...$\,, though the present case of the even $m=2$ in
\ef{m2} will provide us with some surprises, especially for centre
manifold patterns.
 In fact, dealing with
\ef{ex2}, we show how to extend a blow-up solution beyond
singularity towards hence
 guaranteeing existence of a {\em unique} (this is most desirable,
 but not easy to prove) continuation.

\section{
 Self-similar blow-up patterns}
 \label{S2}

\subsection{Blow-up similarity solutions}

For convenience, we reduce the blow-up time to
 $$
 T=0,
 $$
 simply meaning that, by shifting in time, the Cauchy problem for \ef{m2} is considered
 in, say, $(-1,+\iy) \times \re$ with initial data posed at
 $t=-1$.

 Indeed, similarity blow-up is the simplest and most natural
one for scaling invariant equations such as \ef{m2}, where the
behaviour as $t \to 0^-$ is given by a self-similar solution:
 \be
 \label{2.1}
 \tex{
 u_{-}(x,t)=(-t)^{-\frac 1{p-1}} f(y), \quad y= \frac
 x{(-t)^{1/4}},
 }
  \ee
  where  $f(y) \not \equiv 0$ is a  solution of the
  ODE:
  \be
  \label{2.3}
   \tex{
  {\bf A}_-(f) \equiv - f^{(4)}- \frac 14 \, y f' - \frac 1{p-1}\,
  f + |f|^{p-1}f=0 \inB \re, \quad f(\pm\iy)=0.
  }
  \ee

 We recall that, for \ef{m1}, such nontrivial self-similar Type
 I blow-up \ef{2.1},
\be
 \label{2.1S}
 \tex{
 u_{-}(x,t)=(-t)^{-\frac 1{p-1}} f(y), \quad y= \frac
 x{(-t)^{1/2}} \LongA
  f''- \frac 12 \, y f' - \frac 1{p-1}\,
  f + |f|^{p-1}f=0,
 }
  \ee
  is nonexistent. This was first proved in Ad'jutov--Lepin in 1984
  \cite{AL84} (see \cite{GP5} for first applications of the nonexistence to blow-up evolution).
 For \ef{m1N},
   nonexistence  in
  the subcritical Sobolev  range $p \le
 \frac{N+2}{N-2}$ was proved in Giga--Kohn \cite{GK85} in 1985.

  But this is not the case for the RDE--4
 \ef{m2}. Note that \ef{2.3} is a difficult ordinary
 differential
 equation (ODE) with the non-coercive,  non-monotone, and non-potential operators,
 so the problem is
  not variational in any weighted $L^2$-spaces.

In what follows, by $f_0(y)$ we will denote the first monotone
symmetric blow-up profile. We also deal with the second symmetric
profile $f_1(y)$, which seems to be  unstable, or, at least, less
stable than $f_0$.

\begin{figure}
\centering \subfigure[$p=\frac 32, \,\, N=1$]{
\includegraphics[scale=0.52]{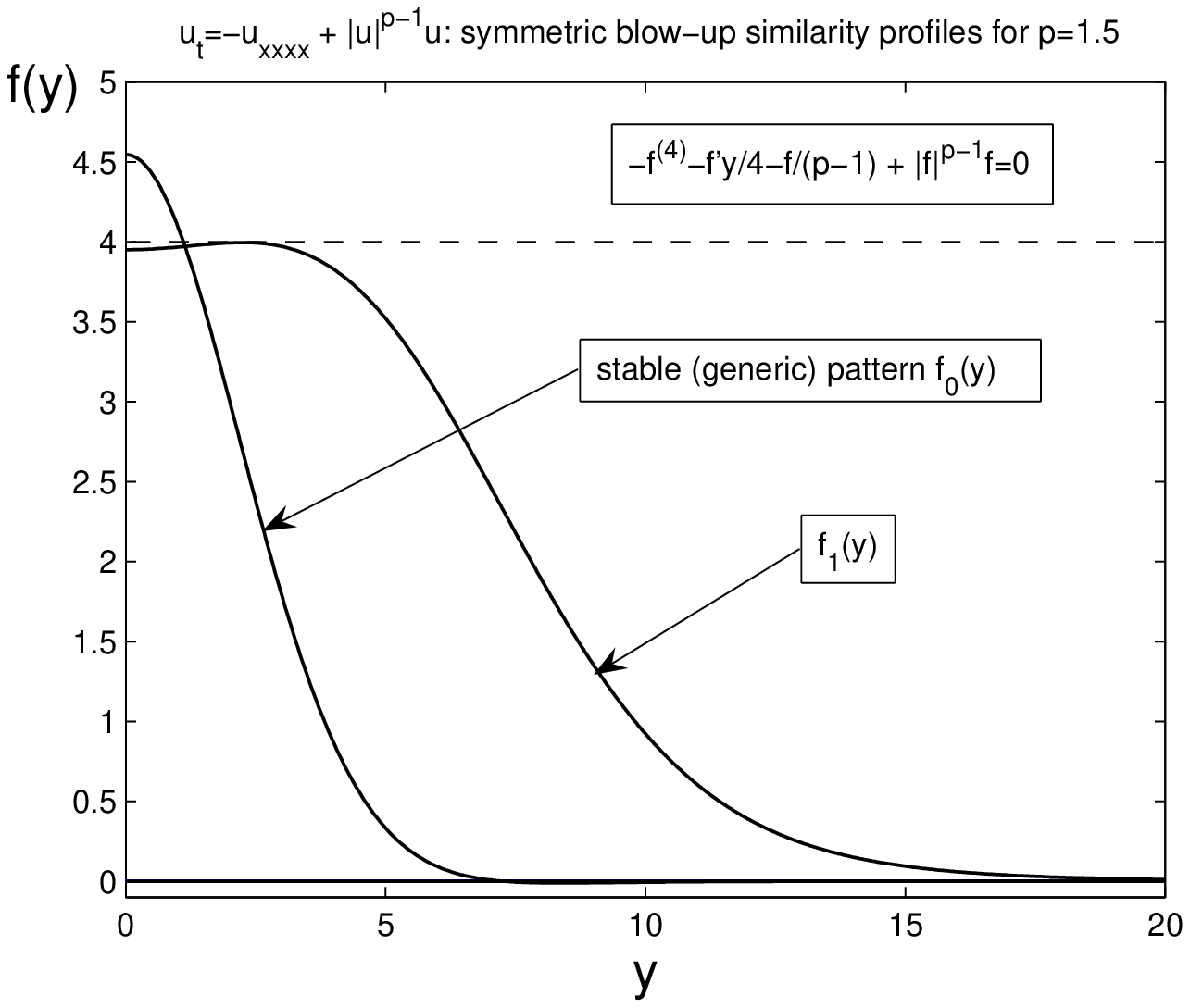}
} \subfigure[$p=2,\,\,N=1$]{
\includegraphics[scale=0.52]{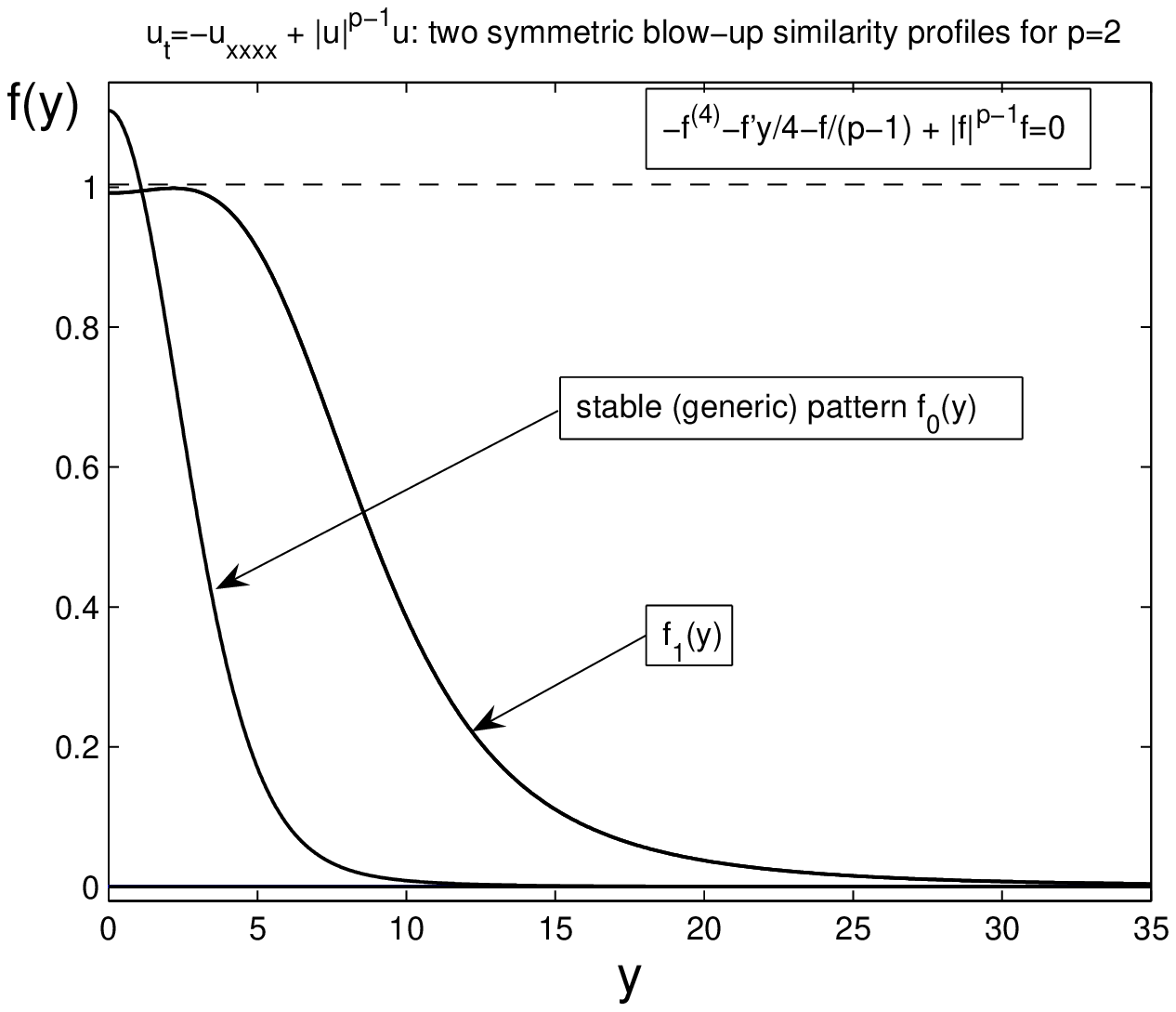}
}
 \vskip -.2cm
\caption{\rm\small Two self-similar blow-up solutions of \ef{2.3}:
$p=1.5$ (a) and $p=2$ (b).}
 \label{F1}
\end{figure}

\begin{figure}
\centering
\includegraphics[scale=0.75]{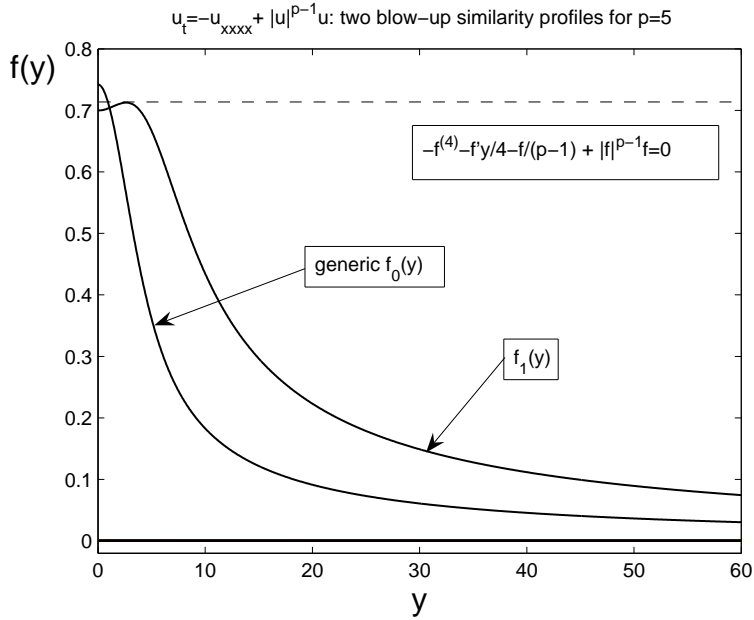} 
\vskip -.3cm \caption{\small  Two self-similar blow-up solutions
of \ef{2.3} for $p=5$.}
 \label{F11}
\end{figure}

The ODE \ef{2.3} was studied in \cite{BGW1} by a number of
analytic-branching and numerical methods. It was shown that
\ef{2.3} admits at least two different blow-up profiles with an
algebraic decay at infinity. We also refer to \cite[\S~3]{GW1} for
further centre manifold-type
 arguments supporting this multiplicity result in a similar
4th-order blow-up problem. Without going into detail of such a
study, we present a few illustrations only (mostly taken from
\cite{GalBlow5}) and will address the essential dependence of
similarity profiles $f(y)$ on $p$. In Figure \ref{F1}, we present
those pairs of solutions of \ef{2.3} for $p=\frac 32$ and $p=2$.
All the profiles are symmetric (even), so satisfy the symmetry
condition
 \be
 \label{2.31}
 f'(0)=f'''(0)=0.
 \ee
No non-symmetric blow-up was detected in numerical experiments
(though there is no proof that such ones are nonexistent: recall
that ``moving plane" and Aleksandrov's Reflection Principle
methods do not apply to \ef{m2} without the MP). Figure \ref{F11}
shows similar two blow-up profiles for $p=5$.


\subsection{Dimension of the ``good" asymptotic bundle at infinity}

 To explain the nature of difficulties in proving existence of
 solutions of \ef{2.3}, let us describe the admissible
 behaviour for $y \gg 1$. It is not difficult to show that there exists a 2D bundle of
 such``good"
 asymptotics (see details in \cite[\S~3.3]{Bl4}): as $y \to +\iy$,
  \be
  \label{dd1}
   \tex{
f(y) = \big[C_1 y^{-\frac 4{p-1}}+... \big]+ \big[C_2 y^{-\frac
23(\frac {p-3}{p-1})}\, {\mathrm e}^{-a_0 y^{ 4/3}}+... \big]
\whereA a_0=3 \cdot 2^{-\frac 83}
 }
 \ee
and $C_1$ and $C_2$ are arbitrary parameters. The same symmetric
bundle exists as $y -\iy$.
 The first square bracket with algebraic decay is connected with
 the dominated linear terms:
  \be
  \label{do9m1}
   \tex{
   - \frac 14 \, f'y- \frac 1{p-1} \, f+...=0 \LongA f(y) = C_1 y^{-\frac 4{p-1}}+...
   \, .
   }
   \ee
 In the second square
bracket in \ef{dd1}, we see a typical WKBJ-type two-scale
asymptotics in ODE theory. The constant $a_0$ in \ef{dd1} is
obtained by substituting into the principal part of \ef{2.3} the
pure exponential term:
 \be
 \label{dd99}
 \tex{
  -f^{(4)}- \frac 14\, y f'+...=0 \withA
f(y) \sim {\mathrm e}^{a y^{ 4/3}} \LongA
 a^3=- \frac 14\, \big( \frac 34\big)^3,
  }
  \ee
  whence the unique (real) root $-a_0$ with a negative real part.
  We do not go into details of the asymptotic expansions  here, since such calculus are well
  known in deriving
 optimal exponential estimates of the fundamental solutions
 of higher-order parabolic equations, which were first
 obtained in Evgrafov--Postnikov (1970) and Tintarev (1982); see
 Barbatis \cite{Barb, Barb04} for key references and updated results.

 Note also that \ef{dd1}  reminds
a typical centre manifold structure of the origin $\{f=0\}$ at $y=
\iy$: the first term in \ef{dd1} is a node bundle with algebraic
decay, while the second one corresponds to ``non-analytic"
exponential bundle around any  of algebraic curves.

Thus, a dimensionally well-posed shooting is characterized as
follows:
 \be
 \label{dd2}
  \fbox{$
 \mbox{{\bf Shooting:} \,\, using {\bf 2} parameters $C_{1,2}$ in (\ref{dd1})
to satisfy {\bf 2} conditions (\ref{2.31}).} $}
 \ee
In view of the analytic dependence of solutions of \ef{2.3} on the
parameters $C_{1,2}$ in the bundle \ef{dd1}\footnote{This is
rather plausible via standard trends of ODE theory \cite{CodL},
but difficult to prove. For positive solutions, it is true, and is
always  straightforward for odd $p=3,5,...$, where the
nonlinearity is analytic.}, the problem \ef{2.3}, \ef{2.31} cannot
have more than a countable set of solutions. Actually, the
numerics and branching-homotopy
 approaches   \cite{BGW1}
 confirm that in wide parameter ranges of $p>1$ and $N \ge
1$, there exist not more than two solutions:
 \be
 \label{dd2N}
 f_0(y) \,\,\, \mbox{with} \,\,\, \{C_{10}(p),C_{20}(p)\},
 \andA f_1(y) \,\,\, \mbox{with} \,\,\, \{C_{11}(p),  C_{21}(p)\}.
  \ee
 More precisely, existence of both  the first (generic blow-up) profile $f_0(y)$ and
  the second one $f_1(y)$
 are
 obtained by a ``$\mu$-bifurcation" approach, when the ODE
 \ef{2.3} is replaced by
  \be
  \label{an1}
   \tex{
   -f^{(4)}- \mu y f'- \frac 1{p-1}\, f + |f|^{p-1}f=0,
    }
    \ee
with a parameter $\mu \ge 0$.
 Linearization about the constant equilibrium and $y$-scaling yield:
  \be
  \label{an2}
   \tex{
   f=(p-1)^{-\frac 1{p-1}} + Y, \,\,\, y=\big(\frac 1{4 \mu}\big)^{\frac
   14}z
    \LongA \BB^* Y + \frac 1{4 \mu}\, Y + O(Y^2)=0,
    }
     \ee
 where $\BB^*=- D^4_z- \frac 14 \, z D_z$ is the well-known Hermite-type operator with
  the discrete spectrum $\s(\BB^*)=\{- \frac k4, \, k=0,1,2,\}$
  and eigenfunctions being generalized Hermite polynomials \ef{psidec},
  \cite{Eg4}; see more
  details in Section \ref{S4.2}.
   It then follows from classic bifurcation theory \cite{VainbergTr} that
  bifurcations occur in \ef{an2} when $-\frac 1{4 \mu}$ gets on the spectrum of $\BB^*$,
  i.e.,
   \be
    \label{an3}
     \tex{
     - \frac 1{4 \mu}=- \frac k 4 \LongA \exists \,\, \mbox{bifurcation points} \,\,
      \mu_k=  \frac 1 k \forA
     k=2,4,6,... \, ,
     }
   \ee
   where by natural symmetry reasons we take into account even
   $k$'s only.
 It turns out that the bifurcations
at $\mu_2= \frac 12$ and $\mu_4= \frac 14$ are responsible for
existence of $f_0$ and $f_1$ respectively; see \cite[\S~5-7]{BGW1}
for details. For convenience, let us note that the first
$\mu$-branch of solutions of \ef{an1} originated at $\mu_2= \frac
12$ (a subcritical bifurcation) is strictly monotone decreasing
for $\mu \in \big(0, \frac 12\big)$ giving the $f_0$ at $\mu=
\frac 14$. The next $\mu$-branch originated at $\mu_4= \frac 14$
(a supercritical bifurcation)  is not monotone and is increasing
on some interval $\mu \in \big(\frac 14, \mu_*\big)$, so there is
another non-zero element on it at $\mu= \frac 14$, which is
precisely the $f_1$. In general, the existence and multiplicity
study of solutions of \ef{2.3} in \cite{BGW1, GalBlow5, GW1} is
difficult and tricky, so here and later on we will need to
essentially rely on careful numerical evidence to check the actual
matching of the flows.

\ssk

The eventual similarity blow-up patterns \ef{2.1} are
characterized by their {\em final  time profiles}: passing to the
limit $t \to 0^-$ in \ef{2.1} and using the expansion \ef{dd1}
yields
 \be
 \label{dd3}
 \mbox{if \,\,$C_1(p) \not = 0$,\,\, then} \quad u_-(x,t) \to C_1
 |x|^{-\frac 4{p-1}} \asA t \to 0^-
  \ee
  uniformly on any compact subset of $\re \setminus \{0\}$.

\subsection{$C_1=0$: final time profile as a measure}

If $C_1 =0$ in \ef{dd1}, i.e., $f(y)$ has  exponential decay at
infinity, then the limit in \ef{dd3} is different and there
appears a measure in the data $u_-(x,0^-)$: in the sense of
distributions,
\be
 \label{dd4}
  \tex{
 C_1(p)  = 0: \quad
  |u_-(x,t)|^{\frac{p-1}4} \to
 E_- \d(x), \,\,
  t \to 0^-; \,\,\, E_-= \int |f|^{\frac{p-1}4} < \iy.
 }
  \ee
Taking into account the sign of the solution, we have a measure in
the ``mass" sense:
\be
 \label{dd4M}
  \tex{
 C_1(p)  = 0: \quad
  |u_-(x,t)|^{\frac{p-5}4}u_-(x,t) \to
 e_- \d(x), \,\,
  t \to 0^-; \,\,\, e_-= \int |f|^{\frac{p-5}4}f \ne 0.
 }
  \ee

\begin{figure}
\centering
\includegraphics[scale=0.75]{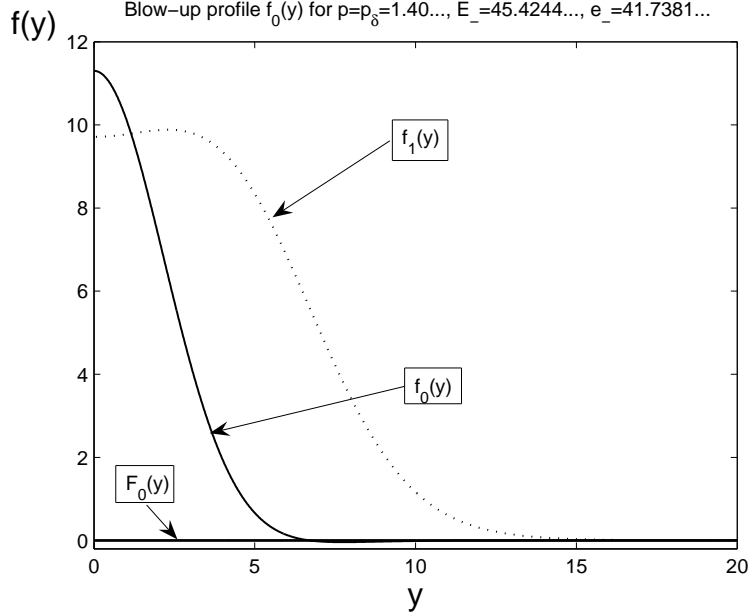} 
\vskip -.3cm \caption{\small  Two self-similar blow-up solutions
of \ef{2.3} for $p=p_\d=1.40...$, when \ef{dd4} holds for $f_0$.}
 \label{FM1}
\end{figure}

 It is  difficult to prove analytically that \ef{dd4} actually takes place at some
 $p=p_\d>1$,
 so numerical methods have been used \cite{GalBlow5} to support this idea.
Namely, the following exponent,  for which \ef{dd4} holds, was
detected:
\be
 \label{rr1}
 p_\d(1)=1.40... \,, \quad \mbox{for which}\quad  E_-=45.4244...\, .
  \, . 
  \ee
 In Figure \ref{FM1}, we show such profile $f_0(y)$ in the case \ef{rr1}, accompanied
 by the second one $f_1(y)$, for which $C_1 \ne 0$.

\subsection{On existence of similarity profiles: classification of blow-up and oscillatory
bundles}

We now prove existence of at least a single blow-up profile $f(y)$
satisfying \ef{2.3}, \ef{2.31}. We perform shooting from $y=+\iy$
 by using the 2D bundle \ef{dd1} to $y=0$, where the symmetry condition \ef{2.31} are posed
 (or to $y=-\iy$, where the same bundle \ef{dd1} with $y \mapsto
 -y$ takes place). By $f=f(y;C_1,C_2)$, we denote the
 corresponding solution defined on some maximal interval
  \be
  \label{ma1}
   \tex{
y \in (y_0, +\iy) \whereA y_0=y_0(C_1,C_2)  \ge -\iy. }
 \ee
If $y_0(C_1,C_2)=-\iy$, then the corresponding solution
$f(y;C_1,C_2)$ is global and can represents a proper blow-up
profile (but not often, see below). Otherwise:
 \be
 \label{ma2}
 y_0(C_1,C_2) >-\iy \LongA f(y;C_1,C_2) \to \iy \asA y \to y_0^+.
  \ee
Note that ``oscillatory blow-up" for the ODE close to $y=y_0^+$:
 \be
 \label{ma3}
  f^{(4)}= |f|^{p-1}f(1+o(1)),
   \ee
where $\lim \sup f(y)=+\iy$ and $\lim \inf f(y)=-\iy$ as $y \to
y_0^+$, is nonexistent. The proof is easy and follows by
multiplying \ef{ma3} by $f'$ and integrating between two extremum
points $(y_1,y_2)$, where  the former one $y_1$ is chosen to be
sufficiently close to the blow-up value $y_0^+$, whence the
contradiction:
 $$
 \tex{
 0< \frac 12\, (f'')^2(y_1) \sim -\frac 1{p+1}\, |f|^{p+1}(y_1)<0.
 }
 $$

We first study this set of blow-up solutions. These results are
well understood for such fourth-order ODEs; see \cite{Gaz06}, so
we omit some details.

\begin{proposition}
 \label{Pr.Inf1}
  The set of blow-up solutions $\ef{ma3}$ is four-dimensional.
   \end{proposition}

\begin{figure}
\centering
\includegraphics[scale=0.65]{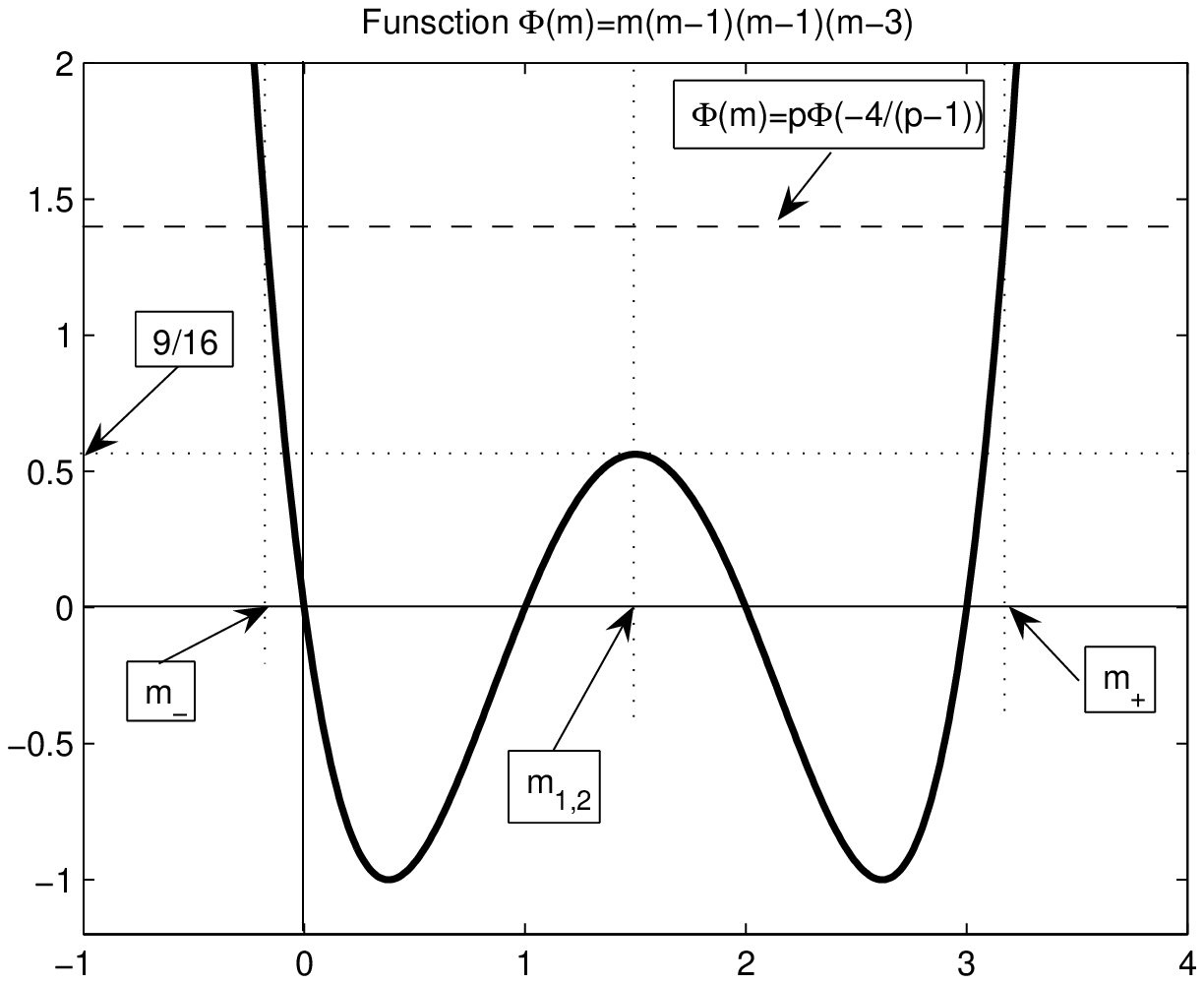} 
\vskip -.3cm \caption{\small The graph of function $\Phi(m)$ in
\ef{ma4}, and towards solutions of \ef{ma5}.}
 \label{FPh}
\end{figure}

\noi{\em Proof.} The first parameter  is $y_0 \in \re$. Other are
obtained  from the principal part of the equation \ef{ma3}
describing blow-up via \ef{2.3} as $y \to y_0^+$.
 We  apply a standard perturbation argument to \ef{ma3}. Omitting
 the $o(1)$-term and assuming that $f>0$, we find its explicit solution
 \be
 \label{ma4}
  \tex{
 f_0(y)=A_0(y-y_0)^{-\frac 4{p-1}}, \,\,\, A_0^{p-1}= \Phi(-
 \frac 4{p-1}), \,\,\, \Phi(m) \equiv m(m-1)(m-2)(m-3).
 }
 \ee
 For convenience, the graph of $\Phi(m)$ is shown in Figure
 \ref{FPh}. Note that it is symmetric relative to $m_0=\frac 32$,
 at which $\Phi(m)$ has a local maximum:
  \be
  \label{loc1}
   \tex{
   \Phi(\frac 32)= \frac 9{16}.
   }
   \ee

By linearization, $f=f_0+Y$, we get Euler's ODE:
 \be
 \label{ma5}
  \tex{
 (y-y_0)^4 Y^{(4)}= p A_0^{p-1} \equiv p \Phi(-
 \frac 4{p-1}).
 }
  \ee
It follows that the general solution is composed from the
polynomial ones with the following characteristic equation:
 \be
 \label{ma6}
  \tex{
  Y(y)= (y-y_0)^m \LongA \Phi(m)=p \Phi(-
 \frac 4{p-1}).
  }
  \ee
 Since the multiplier $p>1$ in the last term in \ef{ma6} and $m=- \frac 4{p-1}$
 is a solution if this ``$p \times$" is omitted, this algebraic equation for $m$ admits a
 unique positive solution $m_+>3$, a negative one $m_- < - \frac 4{p-1}$, which is
 not acceptable by \ef{ma4}, ad two complex roots $m_{1,2}$ with ${\rm Re}\, m_{1,2}= \frac 32>0$.
 Therefore, the general solution of \ef{ma3} about the blow-up one
 \ef{ma4}, for any fixed $y_0$, has a 3D stable manifold. $\qed$

 \ssk

Thus, according to Proposition \ref{Pr.Inf1}, the blow-up
behaviour with a fixed sign  \ef{ma2} (i.e., non-oscillatory) is
generic for the ODE \ef{2.3}. However, this 4D blow-up bundle
together with the 2D bundle of good solutions \ef{dd1} as $y \to
\pm \iy$ are not enough to justify the shooting procedure. Indeed,
by a straightforward dimensional estimate, an extra bundle at
infinity is missing.

\begin{figure}
\centering
\includegraphics[scale=0.65]{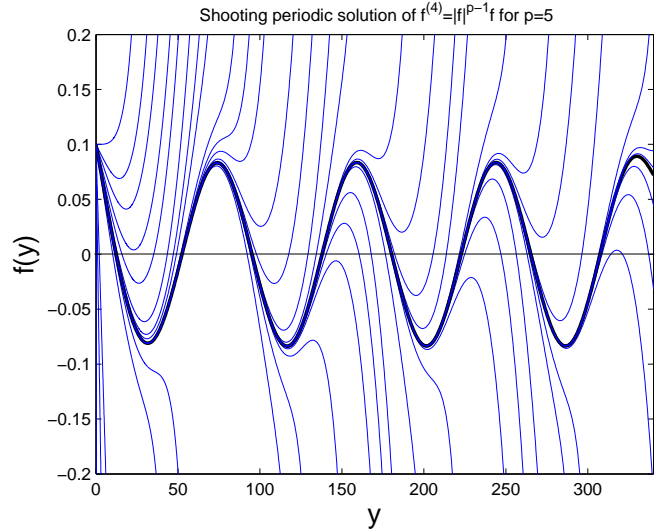} 
\vskip -.3cm \caption{\small A periodic solution of
$f^{(4)}=|f|^{p-1}f$ as a separatrix: $p=5$.}
 \label{FVar}
\end{figure}

To introduce this new oscillatory bundle, we begin with the
simpler ODE \ef{ma3}, without the $o(1)$-term, and present in
Figure \ref{FVar} the results of shooting of a ``separatrix" that
lies between orbits, which blow-up to $\pm \iy$. Obviously, this
separatrix is a periodic solution of this equation with a
potential operator. Such variational problems are known to admit
periodic solutions of arbitrary period.

\begin{figure}
\centering
\includegraphics[scale=0.65]{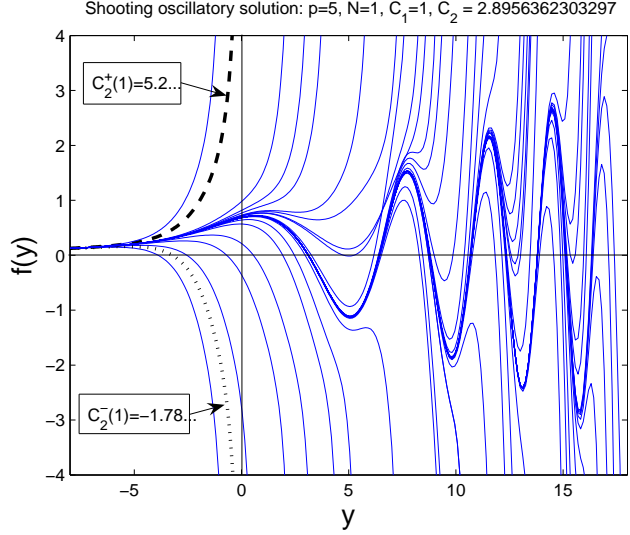} 
\vskip -.3cm \caption{\small Shooting an oscillatory  solution at
infinity  of \ef{2.3}: $p=5$.}
 \label{FF1}
\end{figure}

\begin{figure}
\centering \subfigure[$p=7$]{
\includegraphics[scale=0.52]{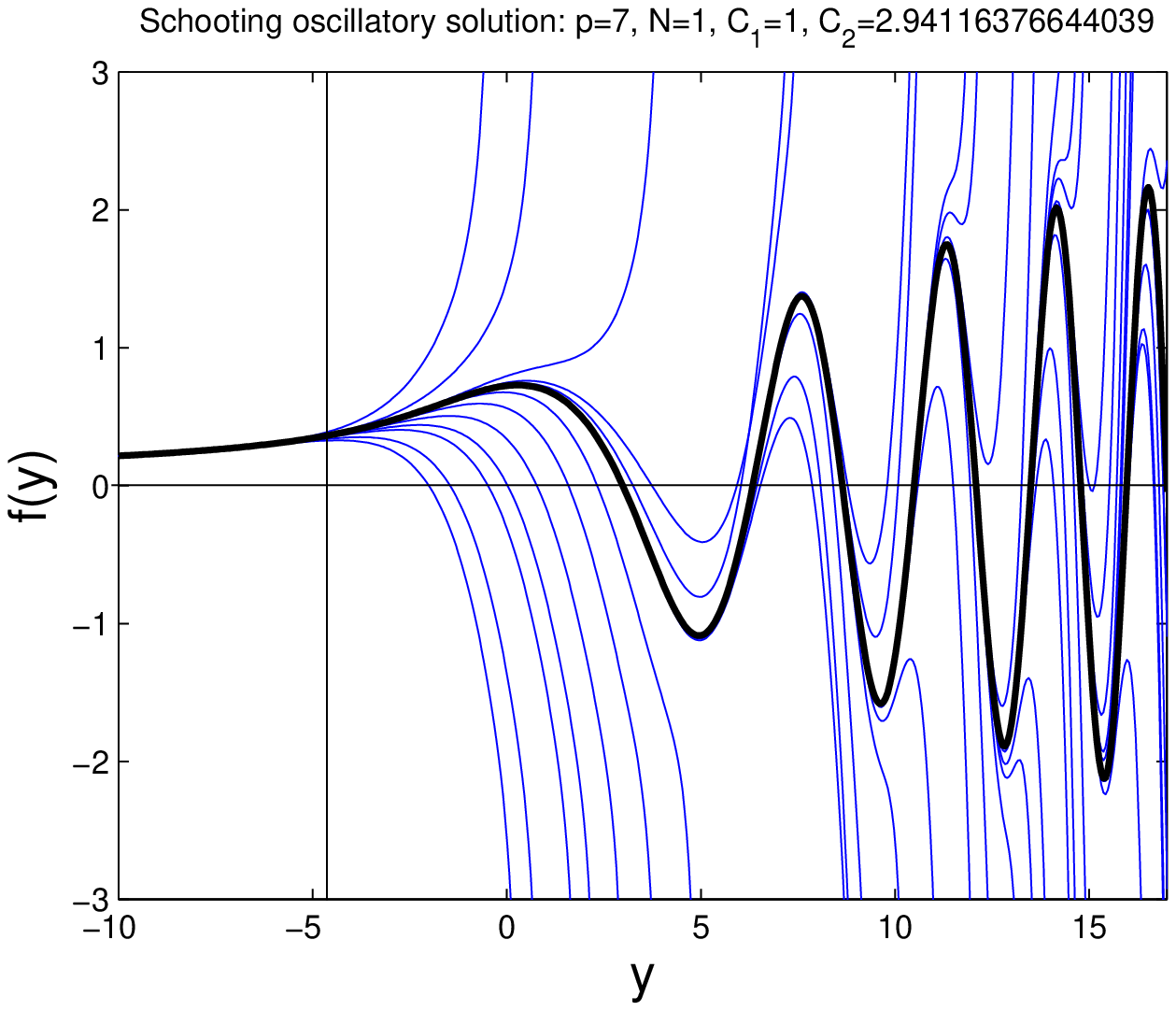}
} \subfigure[$p=2$]{
\includegraphics[scale=0.52]{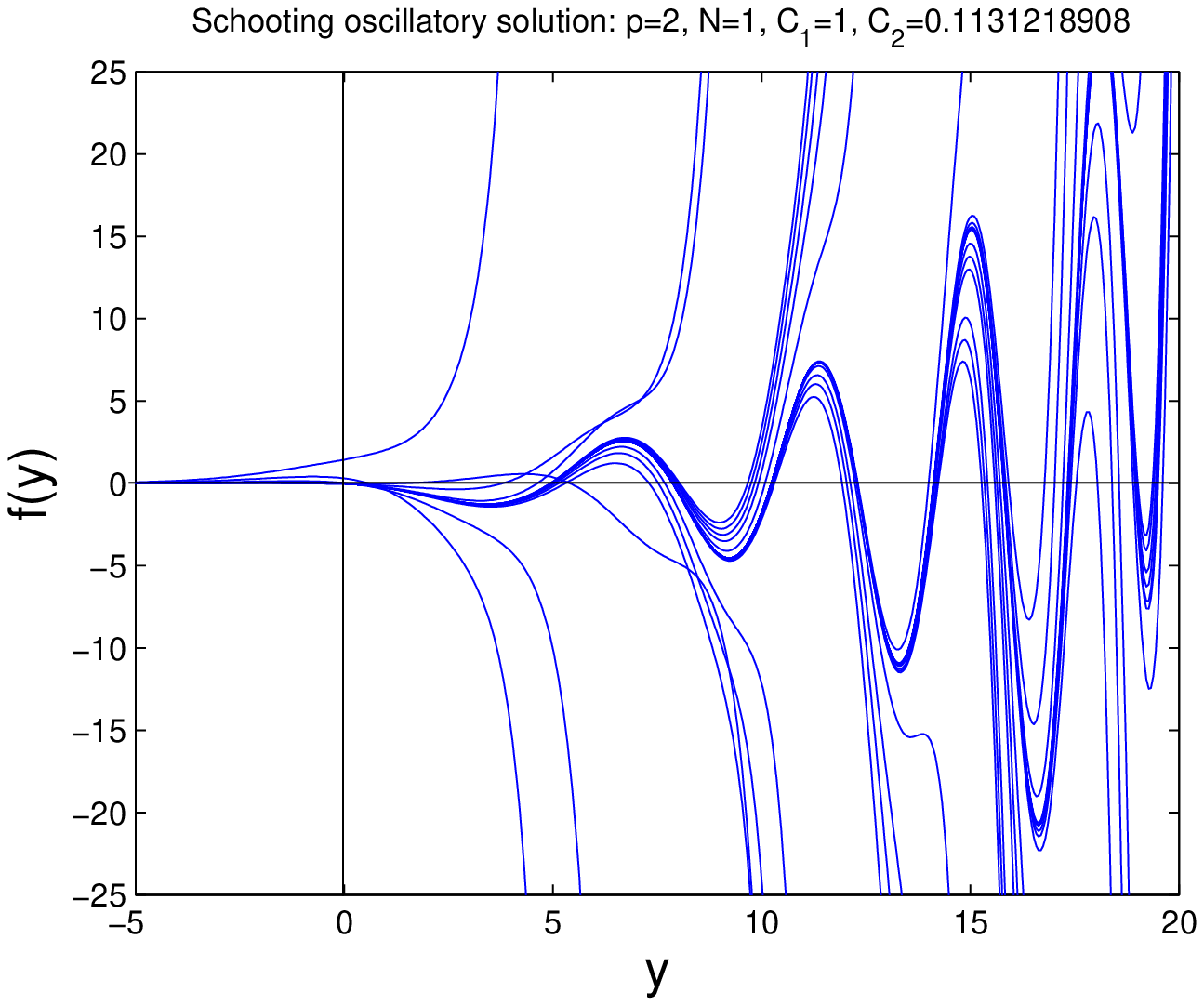}
}
 \vskip -.2cm
\caption{\rm\small  Shooting an oscillatory  solution at infinity
of \ef{2.3}: $p=7$ (a) and $p=2$ (b).}
 \label{FF2}
\end{figure}

Thus, Figure \ref{FVar} fixed a bounded oscillatory (periodic)
solution as $y \to +\iy$. When we return to the original equation
\ef{2.3}, which is not variational, we still are able to detect a
more complicated oscillatory structures at $y=\iy$. Namely, these
are generated by the principal terms in
 \be
 \label{pp1}
  \tex{
 f^{(4)} = - \frac 14\, f'y+ |f|^{p-1}f+... \asA y \to \iy.
 }
  \ee
Similar to Figure \ref{FVar}, in Figure \ref{FF1}, we present the
result of shooting (from $y=-\iy$, which is the same by symmetry)
of such oscillatory solutions of \ef{2.3} for $p=5$. It is easy to
see that such oscillatory solutions have increasing amplitude of
their oscillations as $y \to \iy$, which, as above, is proved by
multiplying \ef{pp1} by $f'$ and integrating over any interval
$y_1,y_2)$ between two extrema. Figure \ref{FF2} shows shooting of
similar oscillatory structures at infinity for $p=7$ (a) and $p=2$
(b). It is not very difficult to prove that the set of such
oscillatory orbits at infinity is 1D and this well corresponds to
the periodic one in Figure \ref{FVar} depending on the single
parameter being its arbitrary period.

 By $C_2^\pm(C_1)$ in Figure \ref{FF1}, we denote the values of
 the second parameters $C_2$ such that, for a fixed $C_1 \in \re$,
 the solutions $f(y;C_1,C_2^\pm)$ blow up to $\pm \iy$ respectively.
These values are necessary for shooting the symmetry conditions
\ef{2.31}.

Thus, overall, using two parameters $C_{1,2}$ in the bundle
\ef{dd1} for $y \gg 1$ leads to a well-posed problem of a 2D--2D
shooting:
 \be
 \label{w1}
  \mbox{find $C_{1,2}$ such that:} \quad
  \left\{
   \begin{matrix}
    y_0(C_1,C_2)=-\iy, \quad \mbox{and} \qquad\qquad\quad\\
    \mbox{no oscillatory behaviour as $y \to -\iy$}.
     \end{matrix}
     \right.
   \ee

Concerning the actual proof of existence via shooting of at least
a single blow-up patterns $f_0(y)$,  by construction and
oscillatory property of the equation \ef{2.3}, we first claim that
in view of continuity relative to the parameters,
 \be
 \label{ccl1}
 \mbox{for any $C_1>0$, there exists $C_2^*(C_1) \in (C_2^-(C_1),C_2^+(C_1))$ such that
 $f'''(0)=0$}.
  \ee
We next  change $C_1$ to prove that at this $C_2^*(C_1)$ the
derivative $f'(0)$ also changes sign. Indeed, one can see that
 \be
 \label{ccl2}
 f'(0;C_1,C_2^*(C_1)) > 0 \forA C_1\ll 1 \andA
f'(0;C_1,C_2^*(C_1)) < 0 \forA C_1\gg 1.
 \ee
 Actually, this means for such essentially different values of
 $C_1$, the solution $f(y;C_1,C_2^*(C_1)$ has first oscillatory ``humps"
 for $y>0$ and $y<0$ respectively. By continuity in $C_1$,
 \ef{ccl2} implies existence of a $C_1^*$ such that
  \be
  \label{ccl3}
  f'(0;C_1^*,C_2^*(C_1^*))=0,
  \ee
  which together with \ef{ccl1} induced the desired solution.
  Overall, the above geometric shooting well corresponds to that applied in the
  standard framework of classic ODE theory, so we do not treat
  this in greater detail. However, we must admit that proving analogously existence of
  the second solution $f_0(y)$ (detected earlier by
  not fully justified arguments of
   homotopy and branching theory and confirmed numerically) is an
   open problem. A more difficult open problem is to show why the
   problem \ef{w1} does not admit non-symmetric (non-even) solutions $f(y)$
   (or does it?).

\section{Self-similar extensions beyond blow-up}
 \label{S3}

\subsection{Global similarity solutions}

Following \cite{GV} devoted to \ef{m1N} and actually using Leray's
blow-up scenario \cite{Ler34}, we suppose that
 the simplest way of extending of self-similar blow-up patterns
for $t>0=T$ is using again similarity global patterns:
 \be
 \label{3.1}
 \tex{
 u_{+}(x,t)=t^{-\frac 1{p-1}} F(y), \quad y= \frac
 x{t^{1/4}},
 }
  \ee
  where  $F$ is a solution of the ODE
  \be
  \label{3.3}
   \tex{
  {\bf A}_+(F) \equiv - F^{(4)}+ \frac 14 \, y F' + \frac 1{p-1}\,
  F + |F|^{p-1}F=0 \inB \re, \quad F(\pm\iy)=0.
  }
  \ee
 In comparison with the blow-up one \ef{2.3}, in \ef{2.3} two linear
terms in the middle have changed their signs. We will show that
this essentially changes the dimension of the asymptotic bundles
and hence overall matching results. Note  that {\em exponentially
decaying} solutions of \ef{3.3}
 were already studied in \cite{GHUni} (these are important for
 extensions of the blow-up of type \ef{dd4}), but now we need other
types of solutions with algebraic decay.

For convenience, we always impose the same symmetry conditions
 \be
 \label{3.31}
 F'(0)=F'''(0)=0,
 \ee
and by $F_0(y)$ we will denote those profiles, which can be
considered for the role of  an extension of the blow-up profiles
$f_0$. Then $F_0(y)$ is not necessarily monotone for $y>0$ and
even its uniqueness is rather questionable; see numerical analysis
below.

\subsection{Dimensional analysis of matching for $C_1 \ne 0$}

We first apply to \ef{3.3} a simple test such as \ef{dd99}
 to get the dimension of the exponentially decaying bundle:
\be
 \label{3.dd99}
 \tex{
F(y) \sim {\mathrm e}^{a y^{4/3}} \LongA
 a^3= \frac 14\, \big( \frac 34\big)^3,
  }
  \ee
 where the algebraic equation for $a \in {\mathbb C}$ admits {\em
 two}
 roots with negative real parts:
 \be
 \label{3a}
  \tex{
 a_\pm = a_0 (-\frac 12 \pm \ii\, \frac{\sqrt 3}2 ) \whereA a_0=3 \cdot 2^{-\frac
 83}.
  }
  \ee
Therefore, instead of \ef{dd1}, the bundle of such asymptotic
orbits is 3D:  as $y \to +\iy$,
  \be
  \label{3dd1}
   \begin{matrix}
F(y) = \big[C_1 y^{-\frac 4{p-1}}+... \big]\ssk\ssk\\
    +\,
\big\{y^{-\frac 23(\frac {p+1}{p-1})}\, {\mathrm e}^{-
\frac{a_0}2\, y^{\frac 43}} \big[ C_2 \cos \big( \frac{a_0 \sqrt
3}2\, y^{\frac 43}\big) + C_3 \sin  \big( \frac{a_0 \sqrt 3}2\,
y^{\frac 43}\big]+... \big\},
 \end{matrix}
 \ee
 where $C_2$ and $C_3$ are arbitrary parameters.
 We consider here the case $C_1 \ne 0$, and will treat the special
 one $C_1=0$, with the data-measure \ef{dd4}, later on.

  The  point is that
the first constant $C_1 \ne 0$ in \ef{3dd1} is fixed by the
blow-up limit \ef{dd3}:
\be
 \label{3dd3}
 u_+(x,0^+) = C_1\,
 |x|^{-\frac 4{p-1}}.
  \ee
 Overall, \ef{dd3} and \ef{3dd3} will provide us with the
 necessary continuity of the unbounded self-similar solution at
the blow-up time in the present case with $C_1 \ne 0$:
 \be
 \label{bl1}
 \fbox{$
 u_-(x,0^-)=u_+(x,0^+)\quad \mbox{for $x \ne 0$, i.e., almost everywhere (a.e.) in $\re$.}
  $}
   \ee
Then $\{f_0,F_0\}$ is called a {\em global extension  similarity
pair}, or simply an {\em extension pair}.

Thus, we again arrive at  a dimensionally well-posed shooting:
 for a fixed values $C_1 \ne 0$,
 \be
 \label{3dd2}
  \fbox{$
 \mbox{{\bf Shooting:} \,\, using {\bf 2} parameters $C_{2,3}$ in (\ref{3dd1})
to satisfy {\bf 2} conditions (\ref{3.31}).} $}
 \ee

\subsection{Numerical analysis of self-similar blow-up extension for $C_1 \ne 0$}

However, the actual solvability of the problem lying behind
\ef{bl1} is  difficult for a rigorous analytic study.
 The main point is that the solvability is very much
 $C_1$-dependent (and mostly is nonexistent for large $C_1$),
 so in many cases we will need again to
rely on careful numerics.

\begin{figure}
\centering
\includegraphics[scale=0.75]{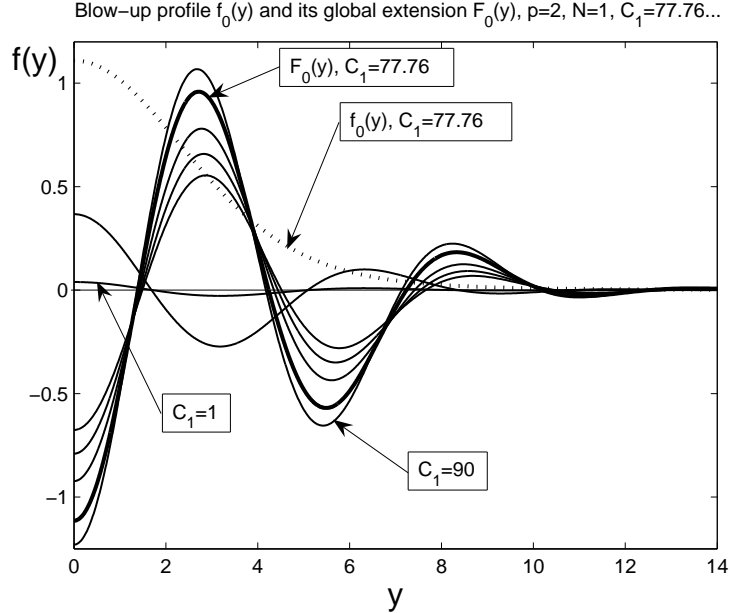} 
\vskip -.3cm \caption{\small  Blow-up profile $f_0(y)$ and the
corresponding global one $F_0(y)$ with the same $C_1=77.76...$ for
$p=2$.}
 \label{F31}
\end{figure}

In Figure \ref{F31}, we show positive results of shooting the pair
$f_0(y)$ and $F_0(y)$ for $p=5$. For comparison, we also put
therein some other global profiles $F(y)$ corresponding to other
values of $C_1$ in \ef{3dd1} from $C_1=1$ up to $C_1=90$, i.e.,
larger than the required 77.76...\,. In the next Figure \ref{F32},
we show the enlarged asymptotic tails of all those profiles.

\begin{figure}
\centering
\includegraphics[scale=0.75]{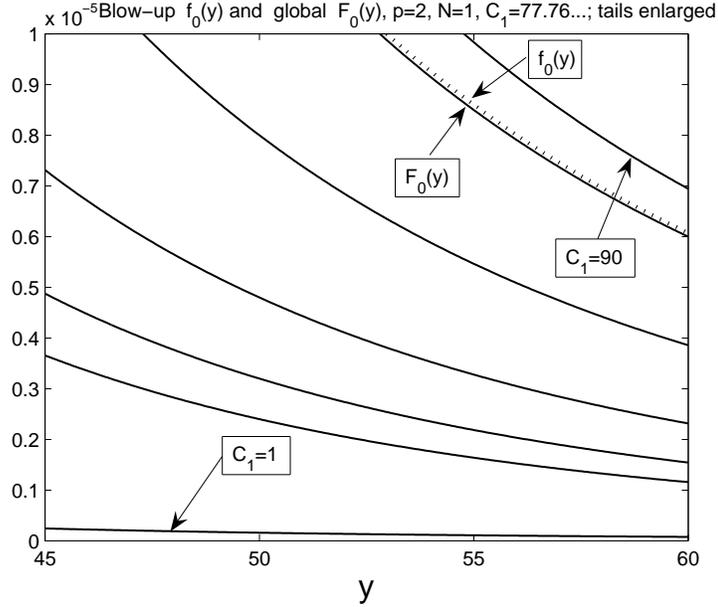} 
\vskip -.3cm \caption{\small  Enlarged tails of blow-up profile
$f_0(y)$ and the corresponding global one $F_0(y)$ from Figure
\ref{F31}  with the same $C_1=77.76...$ for $p=2$.}
 \label{F32}
\end{figure}

\begin{figure}
\centering
\includegraphics[scale=0.75]{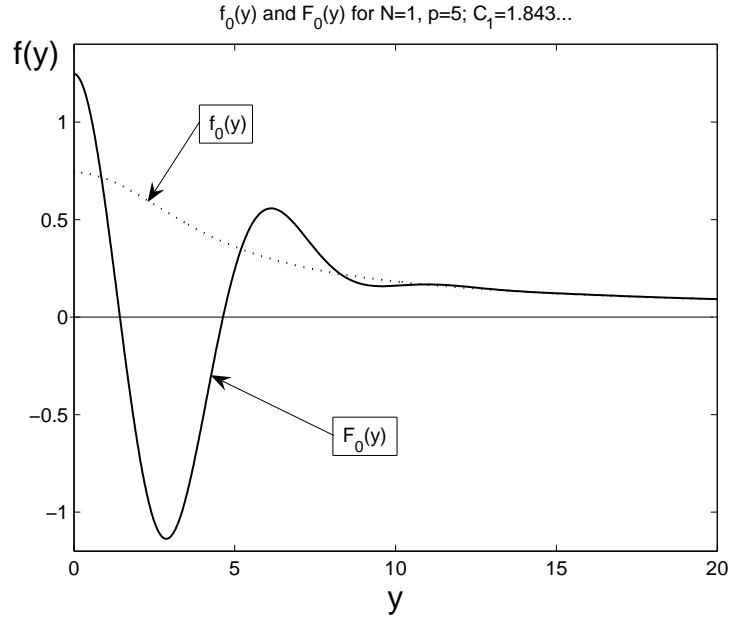} 
\vskip -.3cm \caption{\small  Blow-up profile $f_0(y)$ and the
corresponding global one $F_0(y)$ with the same $C_1=1.843...$ for
$p=5$.}
 \label{F33}
\end{figure}

\begin{figure}
\centering \subfigure[profiles]{
\includegraphics[scale=0.52]{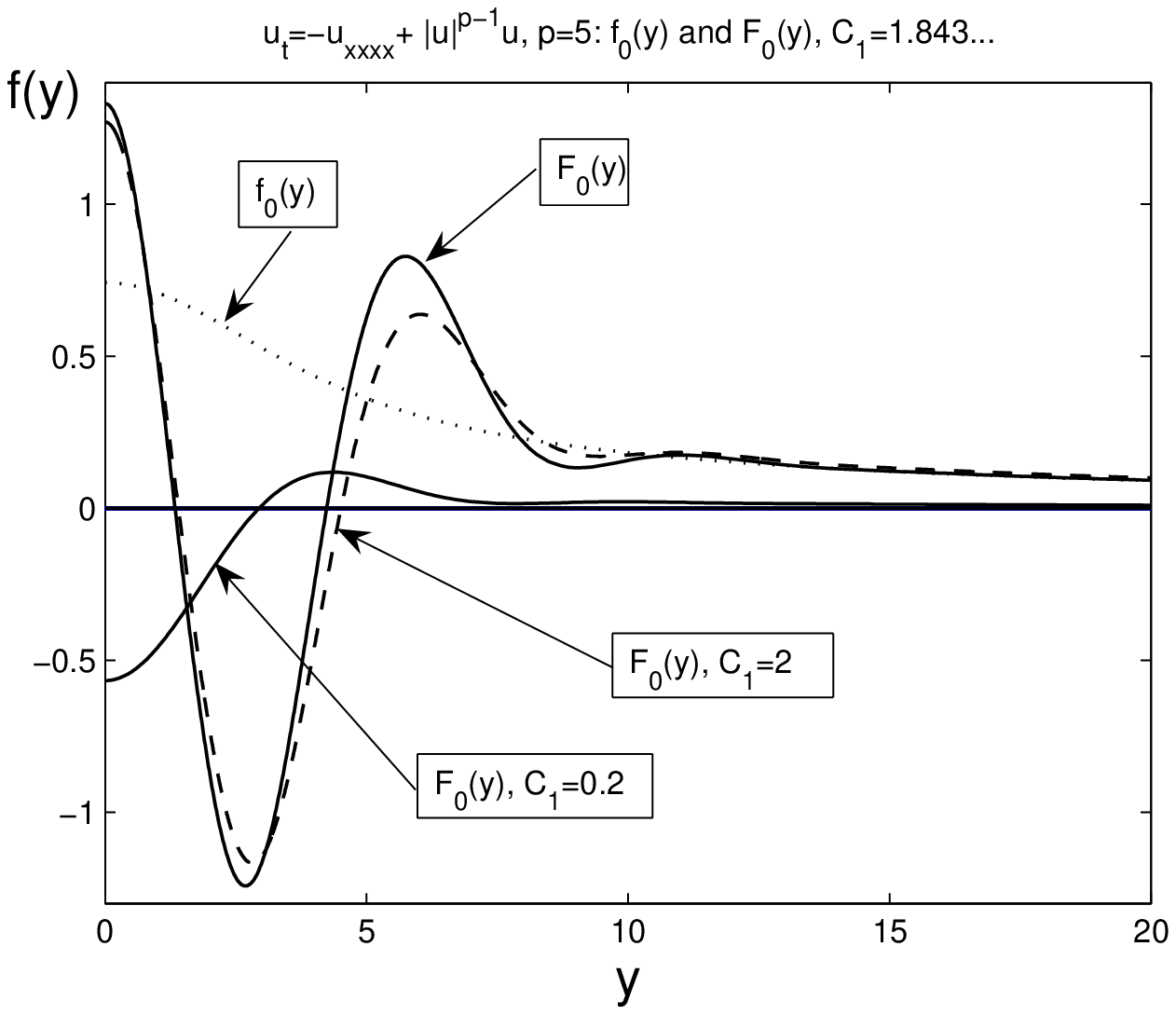}
} \subfigure[enlarged tails]{
\includegraphics[scale=0.52]{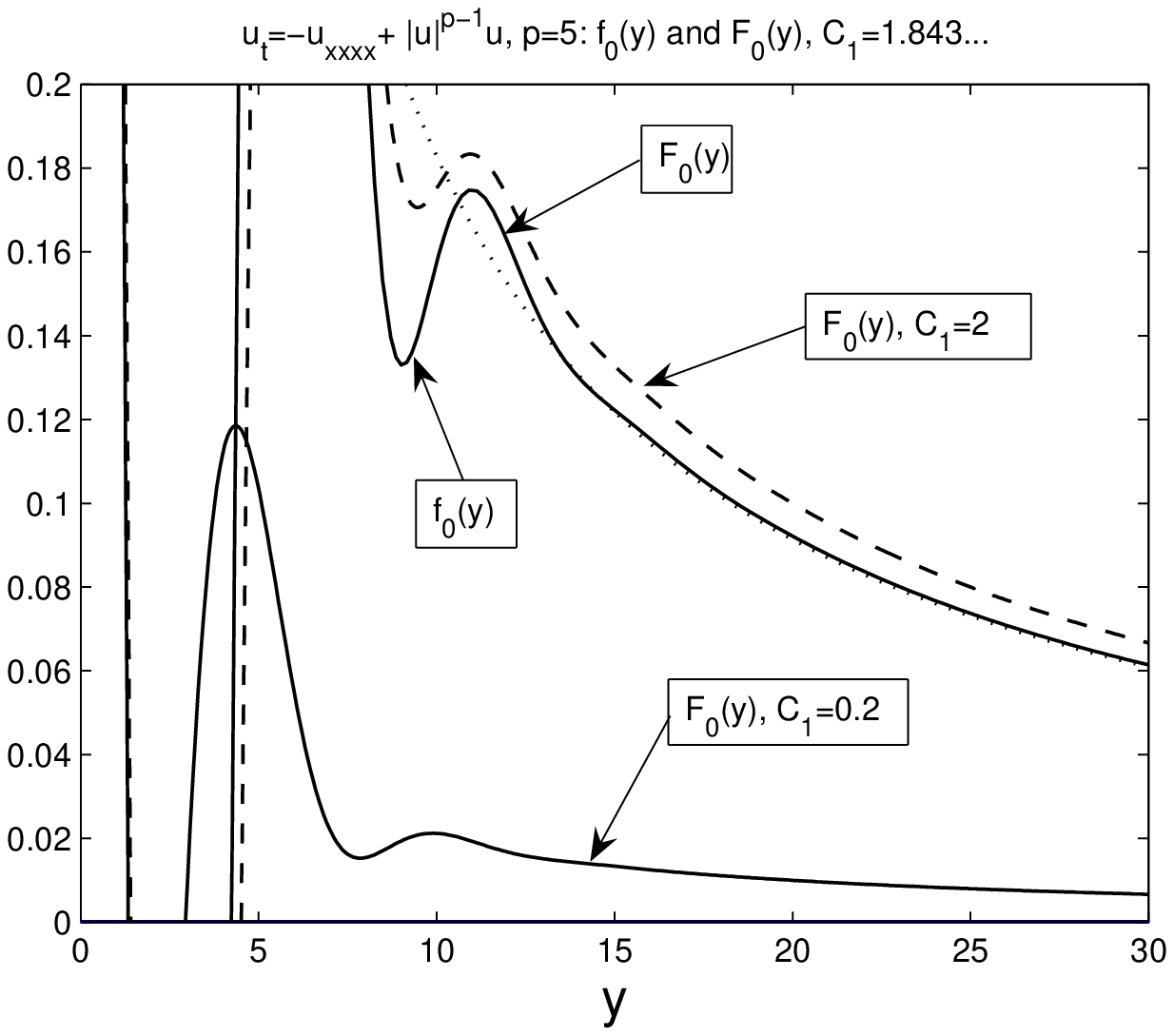}
}
 \vskip -.2cm
\caption{\rm\small  The pair $\{f_0(y),F_0(y)\}$ with
$C_1=1.843...$ for $p=5$ plus extra $F$-profiles.}
 \label{F34}
\end{figure}

Figure \ref{F33} shows the extension pair $\{f_0,F_0\}$ for $p=5$,
where $C_1=1.843...$\,. In Figure \ref{F34}, for the same case of
$p=5$, we present extra global $F$-profiles for smaller $C_1=0.2$
and larger $C_1=2$ values. We should note that convergence for the
global $F$-problem is very slow and we succeeded in getting a few
 reliable numerics only. In all the cases, we have used the  {\tt
bvp4c} solver of the {\tt MatLab} with the enhanced accuracy with
tolerances up to
 \be
 \label{tol1}
 {\rm Tols} \sim 10^{-10}.
   \ee

\begin{figure}
\centering
\includegraphics[scale=0.75]{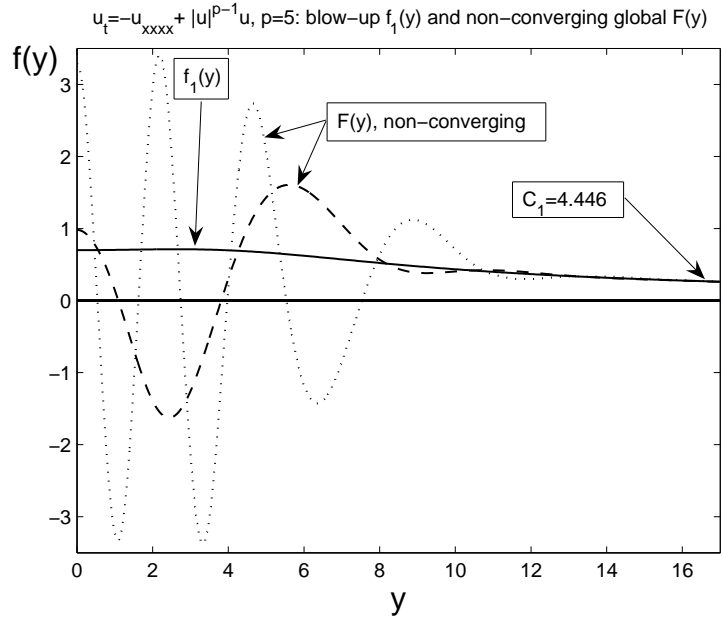} 
\vskip -.3cm \caption{\small Towards nonexistence of $F_1$ for the
second blow-up profile $f_1(y)$, with $C_1=4.446...$ for $p=5$.}
 \label{F35}
\end{figure}

Note that not all the similarity blow-up profiles are assumed to
have a global extension. For instance, Figure \ref{F35}
 explains nonexistence of a global profile $F_1(y)$ for the second
 blow-up one $f_1(y)$ from Figure \ref{F11} for $p=5$. We present
 here the results of non-converging for the equation \ef{3.3} with
 a
 sufficiently large
 $$
 C_1=4.446...\, .
 $$
 Recall that, for $C_1=2$ from Figure \ref{F34}(b), such a profile
 $F_0(y)$ exists, but seems  nonexistent for larger values of $C_1$.
 Possibly, this means that there exists a critical maximal value
 of $C_1$ determining the optimal upper bound for existence of
 $F$-profiles.

\subsection{On extension of a measure for $C_1=0$}

Consider the special case \ef{rr1}, where, as shown in \ef{dd4},
 \be
 \label{mm1}
p=p_\d=1.40... \, : \quad  |u_-(x,0^-)|^{\frac {p-1}4} = E_- \d(x)
\whereA E_-=45.4244...\, .
  \ee
The corresponding blow-up profile $f_0(y)$ has been already shown
in Figure \ref{FM1}.

If $C_1=0$ in \ef{dd1}, one then needs to get the global extension
profile $F_0(y)$  also with $C_1=0$ in \ef{3dd1}. This problem was
studied in \cite{GHUni} by a bifurcation-branching approach. It
was shown that there exists a countable sequence of critical
exponents
 \be
 \label{mm2}
  \tex{
  p_l= 1 + \frac 4{1+l}, \quad l=0,1,2,...\, ,
   }
   \ee
 such that in $\{p>p_l\}$ there exists a global $p$-branch of solutions
 $F_l(y)$ of \ef{3.3}, that can be extended unboundently as $p \to +\iy$.
The behaviour of the branches  near bifurcation points \ef{mm2}
was proved to be:
 \be
 \label{mm3}
  \tex{
  F_{0l}(y)= \g_l (p-p_l)^{\frac {1+l}4} \, \big[\psi_l(y)+o(1)\big]
  \asA p \to p_l^+; \quad l=0,1,2,...\, ,
  }
  \ee
  where $\g_l>0$ are some  constants and $\psi_l(y)$ are eigenfunctions
  \ef{gen11} of the  operator
  \ef{BBB}.

Before constructing the extension pair for the case \ef{mm1}, note
that this $p_\d$, within our accuracy, is close to the spectrum
\ef{mm2}, and, precisely,
 \be
 \label{mm4}
  \tex{
 p_9= 1+ \frac 4{10}=1.40 .
  }
  \ee
Note also that the critical exponents \ef{mm2} are concentrated
about $p=1^+$ for large $l$, and, in view of the asymptotics
\ef{mm3}, the $p$-branches are rather plain close to $p=p_\d$.
Overall, this shows that there exist many global profiles
$F_{0l}(y)$ satisfying \ef{mm1}.

\begin{figure}
\centering
\includegraphics[scale=0.75]{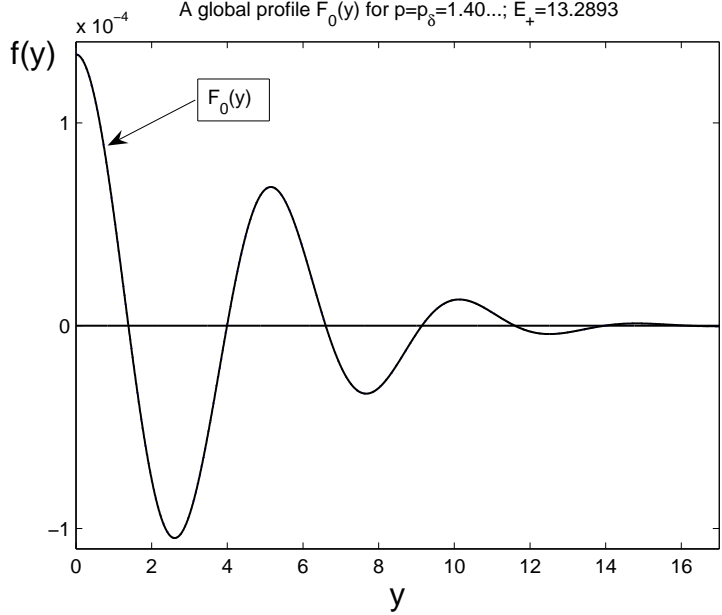} 
\vskip -.3cm \caption{\small A typical global similarity profile
$F_0(y)$ for $p=p_\d=1.40...\,$.}
 \label{FM2}
\end{figure}

 Figure \ref{FM2} shows
a typical global similarity profile $F_0(y)$, though by the
branching \ef{mm3}, there exists a countable set of such solutions
of \ef{3.3} with exponential decay. Therefore, it is difficult to
identify which $p$-branch this profile $F_0$ belongs to. The
measure characteristic of this $F_0(y)$ in Figure \ref{FM2} is:
 \be
 \label{mm5}
 E_+=13.2893... \, .
  \ee
  The computations have been performed with the enhanced tolerances
  \ef{tol1}, since many profiles $F_0$ are small enough due to \ef{mm3}
  for $p \approx p_\d$. On the other hand, by the same branching,
  for $l \gg 1$ with $p_l$ more closer to $1^+$, there exist
 $F_0$'s with arbitrarily large $E_{+l}$.


Indeed, we see that the value
\ef{mm5} for a global extension does not match the blow-up one in
\ef{mm1}:
 \be
 \label{mm6}
 p=p_\d=1.40...\, : \quad E_- > E_{+}.
  \ee
Recall that, due to \ef{mm2}, there also exist global profiles
$F_{0l}$ such that $E_- < E_{+l}$. However, all these
discrepancies do not undermine a possibility of extension of the
data \ef{dd4} in a self-similar way given by \ef{3.1}. An
effective scenario of such a matching was obtained and confirmed
in \cite{EGW1} for the related {limit unstable
 Cahn--Hilliard equation} \ef{mm7},
   which preserves the mass of the $L^1$-solutions.
 Namely, it was shown \cite[\S~4]{EGW1} that the rescaled blow-up solution
 can converge in the local topology (\ref{un1}) to a self-similar
 profile with a different mass ($L^1$-norm for positive solutions). In practice, this
 means that close to blow-up time (as $t \to 0^-$), negative humps
 are created that eventually disappear at infinity, and hence do not
 violate the matching.
 This ``mass
drift" is performed in the rescaled variables as $t \to 0^-$, so
does not mean existence of the actual mass transferring mechanism
in the $x$-variable that might be unavailable in the PDE. There is
no a rigorous proof.

We expect that a similar mechanism can be explored to neglect the
mass ``defect", and, for $t>0$, we may observe a global profile
$F_{0l}(y)$, which  minimizes the loss of the mass $E_- - E_{+l}$.
The positive mass defect will then  create
 two negative/positive humps that disappear at infinity in the rescaled variable $y$ as $t \to 0^+$.
 In this mass interpretation, we have to use the measure-like
 data obtained in the limit \ef{dd4M}.
 Recall that the whole
set of global profiles $\{F_{0l}\}$ at $p=p_\d$ is countable, so
that this minimization makes sense.
 It looks then like a certain ``mass-discontinuity" is assumed at
 $t=0^+$, though this happens in the $y=x/t^{1/4}$ variable and
 hence the real discontinuity is not available in the $x$-one,
 which is prohibited for this parabolic flow.

 There are several open problems here, and the
further analysis of the actual extension pair $\{f_0,F_{0l}\}$
based on enhanced {\em PDE numerics} is quite necessary and
desirable to clarify this difficult transition blow-up singularity
phenomenon.


\subsection{Towards extended semigroup theory}

Thus, the above results, at least at a qualitative formal level,
imply that self-similar blow-up for \ef{m2} (which is assumed to
be generic and most structurally stable; see the next section)
 admits a proper extension beyond blow-up time. Moreover, we also can claim
 that, for $C_1 \ne 0$:
  \be
  \label{cl1}
   \mbox{there exists at most a finite number of extension pairs $\{f_0,F_0\}$.}
 \ee
 For $C_1=0$, i.e., for $p=p_\d$, due to the critical bifurcation
 exponents \ef{mm2}, we expect a countable sequences of the
 extension pairs, but the discrepancy of their ``masses" at $t=0$
 could reinforce extra evolution mechanism of choosing the right
 ones. This  remains an open problem.

 Note that, in general,  we cannot guarantee that a pair is unique.
 Nevertheless, in view of \ef{cl1}, even in the case of finite
 multiplicity of the  profiles $F_0$, there is a hope
 that
  \be
  \label{cl2}
   \mbox{there exists a ``minimal" extension pair $\{f_0,\bar F_0\}$,}
 \ee
 with $\bar F_0$ being ``ordered" in a certain geometric-metric
 sense,
 e.g., $F_0$ with the same  $C_1$ but with a most
 ``less oscillatory" structure. Recall that $f_0$ is always
 monotone.

This {\em minimality property} can play a key role in an  attempt
(seems rather naive) to constructing of extended semigroup theory
for blow-up solutions of \ef{m2} by using standard ideas of
parabolic regularizations via smooth solutions $\{u_\e\}$ (cf.
\cite{GV} for \ef{m1}),
 \be
 \label{reg1}
  \tex{
  u_\e: u_t=-u_{xxxx} + \var_\e(u) \whereA \var_\e(u)=
  \frac{|u|^{p-1}u}{1+\e |u|^{p-1}u} \to |u|^{p-1}u \asA \e \to
  0^+
   }
 \ee
 uniformly on compact subsets. Since $\var_\e(u)$ is globally Lipschitz continuous,
 the CP for \ef{reg1} with the same
 data $u_0$ has a unique global classical solution $u_\e(x,t)$.
 Then a proper (minimal) solution of \ef{m2} can be formally defined as a
 hypothetical  limit
 \be
 \label{reg2}
  \tex{
  \bar u(x,t) = \lim_{\e \to 0} u_\e(x,t) \inB \re \times
  (-1,\iy),
   }
   \ee
   where the limit is understood in a pointwise sense, since
   $u(x,t)$ is unbounded for $t \ge 0$. Moreover,  no more restrictions on the
   topology of convergence in \ef{reg2} can be imposed, since, for
   general solutions, the {\em blow-up set}:
    \be
    \label{se1}
    B[u_0](t)=\{x \in \re: \quad |u(x,t)|=\iy\},
     \ee
    has an unknown structure for $t \ge T$ (hopefully, of zero measure, which is
    not proved yet). Recall that even for much simpler nonlinear
    heat equations including \ef{m1N}, the limit extended
    semigroups are discontinuous in time in general; see
    \cite[Ch.~7]{GalGeom} for examples.

   Proving existence of the limit \ef{reg2} (even along a subsequence) and
   checking how it is related to the minimal $\bar F_0$-extension
   of the blow-up self-similar solution \ef{2.1} are difficult
   open problems (possibly, non-solvable in general). Anyway, it is indeed surprising that there is
   still a slight hope that extended semigroup theory of blow-up
   solutions for \ef{m2} can be at least partially developed along
   the lines of that for the second-order parabolic PDEs such as  \ef{m1} \cite{GV, GalGeom},
   where the Maximum Principle was always key.

\subsection{On generalizations to $\ren$}

\begin{figure}
\centering \subfigure[$N=3$, $C_1=2.328...$]{
\includegraphics[scale=0.52]{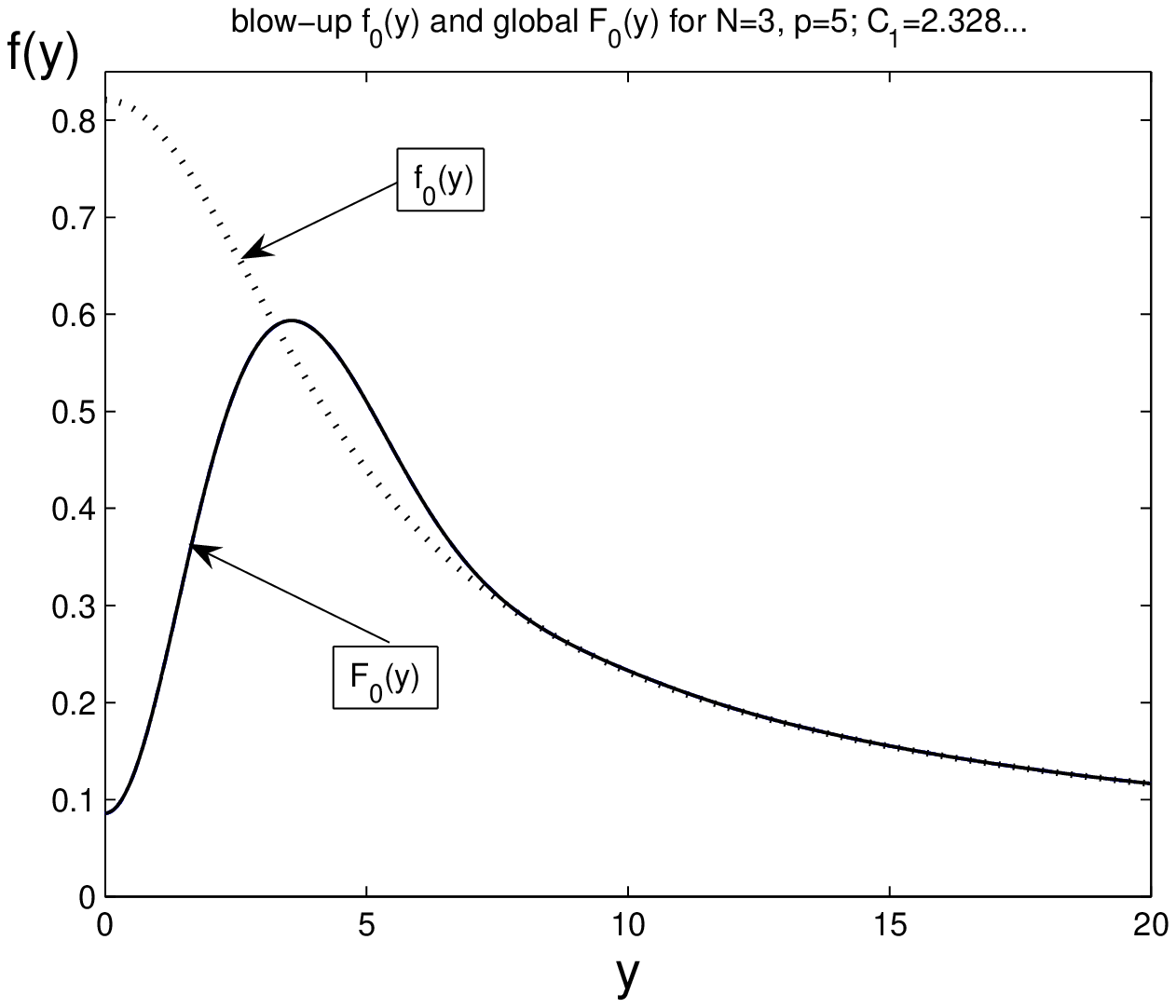}
} \subfigure[$N=6$, $C_1=2.91...$]{
\includegraphics[scale=0.52]{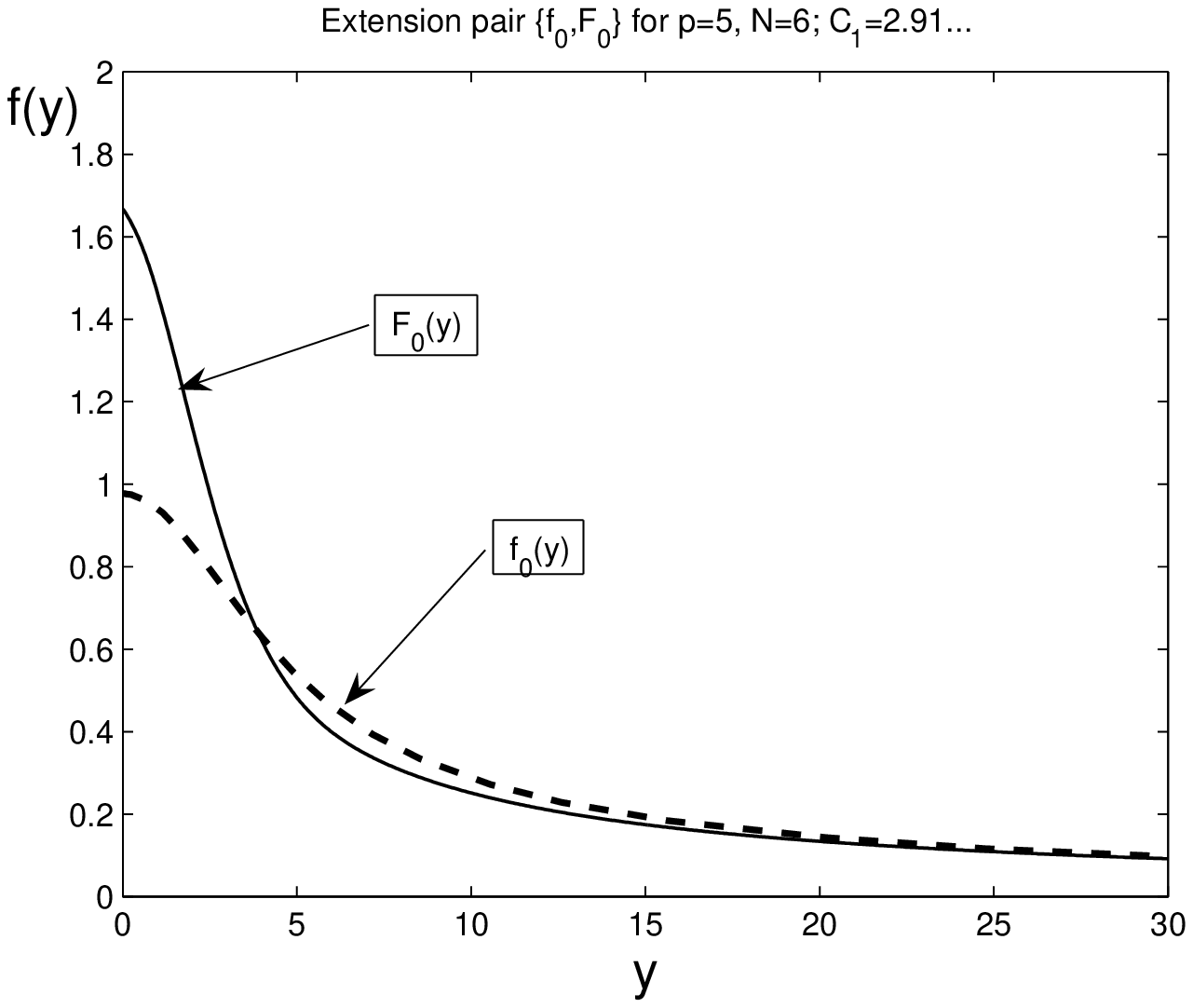}
}
 \vskip -.2cm
\caption{\rm\small  The pair $\{f_0(y),F_0(y)\}$ for $p=5$ and
$N=3$ (a) and $N=6$ (b).}
 \label{F36}
\end{figure}

Finally, as a comment and an introduction to a future research, in
Figure \ref{F36} we show the extension pair $\{f_0,F_0\}$ for the
equation in $\ren$,
 \be
 \label{r3}
 u_t=-\D^2 u + |u|^{p-1}u,
  \ee
in the case $p=5$ for $N=3$ (a) and 6 (b). Note the clear
difference in
 the geometry of the corresponding global
$F_0$-profiles in (a) and (b). This is again a sign showing that
the extension profiles $F_0$ can be nonunique.

The similarity solutions $u_\pm$ remain the same and the ODEs for
radial patterns $\{f_0,F_0\}$ are easily obtained; see
\cite{GalBlow5} for details on blow-up ones $u_-(x,t)$. The ODE
problems for $f_0$ and $F_0$ in the radial geometry are as follows
(now $y$ stands for $|y|>0$):
 $$
  \left\{
 \begin{matrix}
-f^{(4)} - \frac {2(N-1)}y\, f'''- \frac
 {(N-1)(N-3)}{y^2}\, f'' + \frac{(N-1)(N-3)}{y^3}\, f' - \frac
 14\, y f'
 - \frac 1{p-1}\,
  f + |f|^{p-1}f=0, \ssk\\
  f'(0)=f'''(0)=0, \quad f(y)= C_1y^{-\frac4{p-1}}+... \asA y
  \to \iy; \qquad\qquad\qquad\qquad\quad
   \end{matrix}
    \right.
  $$ 
 $$
\left\{
 \begin{matrix}
-F^{(4)} - \frac {2(N-1)}y\, F'''- \frac
 {(N-1)(N-3)}{y^2}\, F'' + \frac{(N-1)(N-3)}{y^3}\, F' + \frac
 14\, y F'
 + \frac 1{p-1}\,
  F + |F|^{p-1}F=0, \ssk\\
  F'(0)=F'''(0)=0, \quad f(y)= C_1 y^{-\frac4{p-1}}+... \asA y
  \to \iy. \qquad\qquad\qquad\qquad\qquad
   \end{matrix}
    \right.
   $$ 
 The asymptotic bundles \ef{dd1} and \ef{3dd1} remain analogous to the case $N=1$.

\begin{figure}
\centering \subfigure[$N=8$, $C_1=2.22...$]{
\includegraphics[scale=0.52]{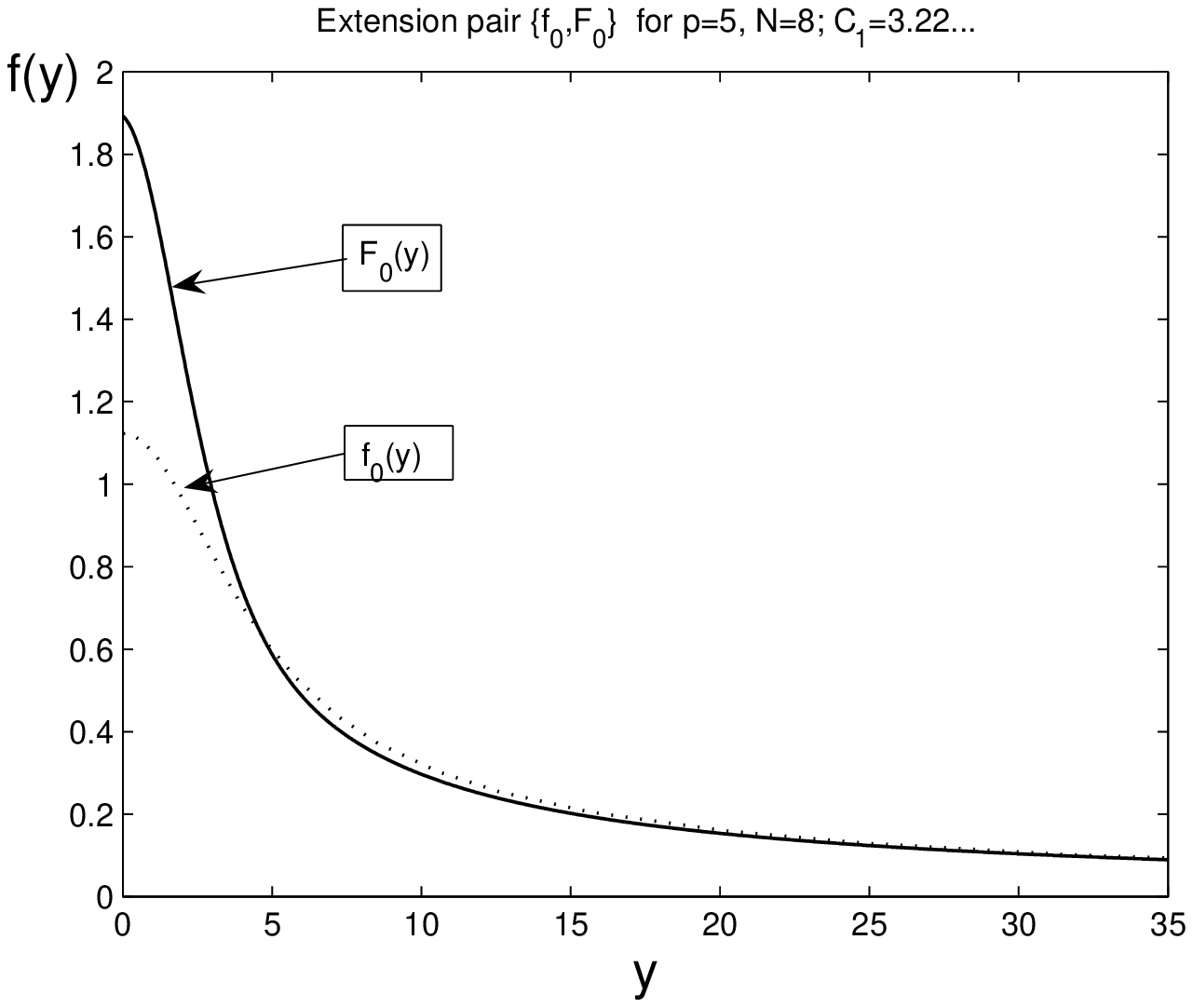}
} \subfigure[$N=9$, $C_1=3.38...$]{
\includegraphics[scale=0.52]{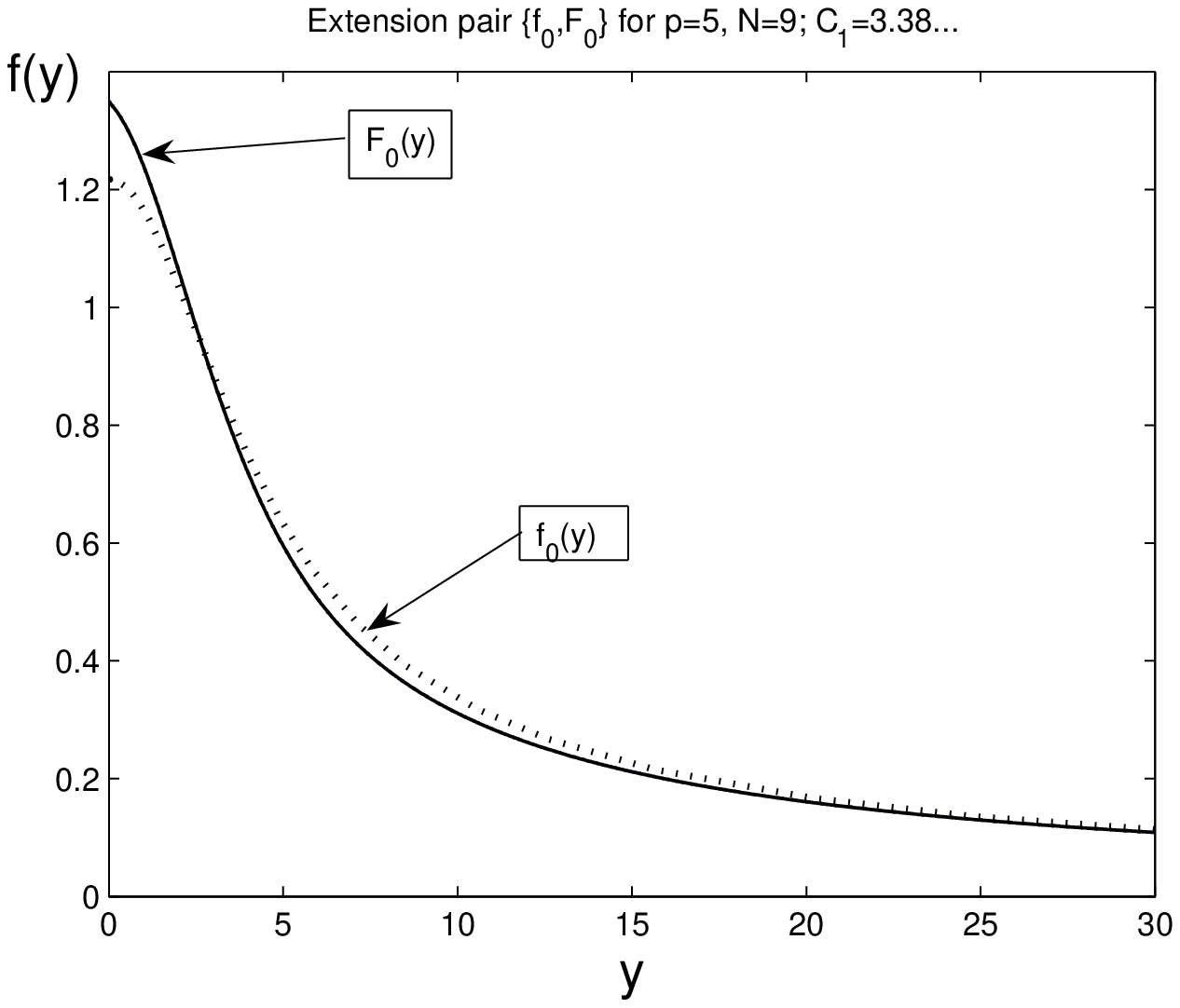}
}
 \vskip -.2cm
\caption{\rm\small  The pair $\{f_0(y),F_0(y)\}$ for $p=5$ and
$N=8$ (a) and $N=9$ (b).}
 \label{F37}
\end{figure}

\begin{figure}
\centering
\includegraphics[scale=0.75]{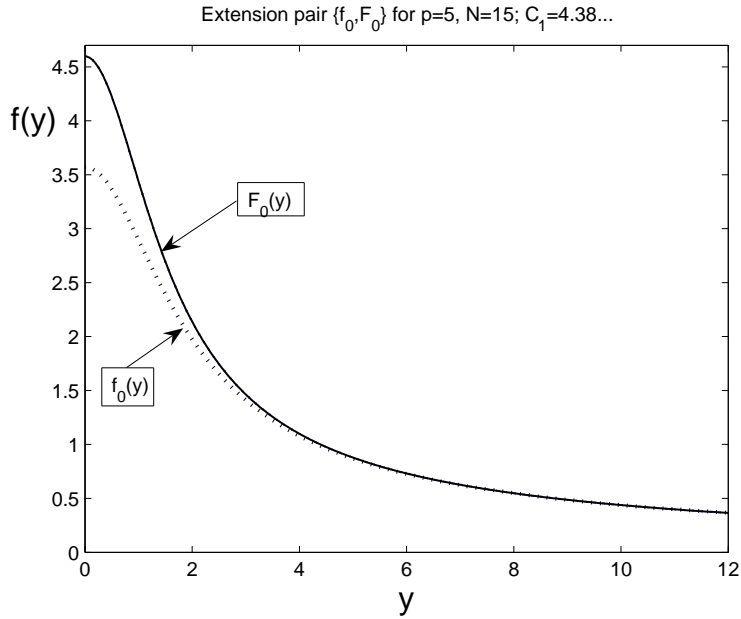} 
\vskip -.3cm \caption{\small  The pair $\{f_0,F_0\}$ for $p=5$,
$N=15$.}
 \label{F38}
\end{figure}

In Figure \ref{F37}, we show the pairs $\{f_0,F_0\}$ again for
$p=5$ in dimensions $N=8$ (a)  and $N=9$ (b).
 Finally, in Figure \ref{F38}, the extension pair is shown for
 $p=5$ and $N=15$. Analogously to $N=1$, the convergence
 for the global profile $F_0(y)$ is very slow, possibly in view of
 the multi-dimension of the bundle \ef{3dd1}. However, for such profiles, we always observed
 a multiple convergence to such profiles from various initial data
 and/or varied accuracy and tolerances used in the {\tt bvp4c} solver.

\section{Non self-similar ``linearized" patterns:
complete blow-up}
 \label{S4}

A countable set of such  non-self-similar blow-up patterns
 for \ef{m2} were formally constructed in \cite{Gal2m}; see also
 \cite{GalBlow5} for extra details and necessary historical comments.
 Below, we briefly describe some necessary new features of these
 patterns to reveal the main reasons for their complete blow-up.
For the RDE--4, there is no  hope to get a full rigorous
justification of existence of such non-self-similar blow-up
scenarios, so we feel free to perform a detailed formal
construction.
 More advanced expansion and matching techniques
 in this direction can be found in
\cite{Gal2m} and also
 \cite{GalCr}, where such a construction applied to
non-singular absorption phenomena (regular flows with
 no blow-up), so a full mathematical justification is available therein.

\subsection{Nonstationary rescaling}

 Dealing with non-self-similar blow-up,
instead of \ef{2.1}, we use the full similarity scaling:
 \be
 \label{3.1MM}
  \tex{
 u(x,t)=(-t)^{-\frac 1{p-1}}v(y,\t), \quad y= \frac
 x{(-t)^{1/4}}, \quad \t = - \ln(-t) \to +\iy \,\, \mbox{as}\,\, t \to 0^-.
}
  \ee
 Then $v(y,\t)$ solves the following parabolic equation:
  \be
  \label{3.2}
   \tex{
  v_\t= {\bf A}_-(v) \equiv - v_{yyyy} - \frac 14\, y v_y -\frac 1{p-1}\, v +
  |v|^{p-1}v\inB \re \times \re_+,
  }
   \ee
   where ${\bf A}_-$ is the stationary  operator in
   \ef{2.3}, so that similarity profiles are simply
   stationary solutions of \ef{3.2}.

\subsection{Linearization and spectral properties}
 \label{S4.2}

Performing the standard linearization about the constant
equilibrium in the equation \ef{3.2} yields a perturbed equation:
 \be
 \label{4.1}
 \tex{
 v=f_*+Y, \,\, f_*=(p-1)^{-\frac 1{p-1}} \LongA
Y_\t= (\BB^*+I) Y  +\DD(Y),
 }
 \ee
where $\DD(Y)=c_0 Y^2+...$, $c_0=\frac {p}2(p-1)^{\frac{1}{p-1}}$,
is a quadratic perturbation as $Y \to 0$ and
  \be
  \label{4.2}
  \tex{
  \BB^*=-D_y^4 - \frac 1{4}\, y D_y \inB L^2_{\rho^*}(\re),
  \,\,\,\rho^*(y)={\mathrm e}^{-a|y|^{4/3}}, \quad a \in \big(0,3 \cdot
 2^{- \frac {8}3}\big),
  }
  \ee
  is the {\em adjoint Hermite operator} with some good spectral
  properties \cite{Eg4}:

  \begin{lemma}
\label{lemSpec2} 
  $\BB^*: H^{4}_{\rho^*}(\re) \to
L^2_{\rho^*}(\re)$ is a bounded linear operator with the spectrum
 \be
\label{SpecN}
 \tex{
 \s(\BB^*) = \{ \l_k =  -\frac{k}4, \,\,\,
k=0,1,2,...\} \quad \big(= \s(\BB), \,\, \BB=-D_y^4 + \frac 1{4}\,
y D_y+ \frac 14 \, I \big). }
   \ee
     Eigenfunctions
$\psi_k^*(y)$ are $k$th-order generalized Hermite polynomials:
\be
\label{psidec}
 \tex{
\psi_k^*(y) =  \frac {1}{ \sqrt{k !}}\big[ y^k +
\sum_{j=1}^{[k/4]} \frac 1{j !}(D_y)^{2 j} y^k \big],
 \quad k=0,1,2,...\, ,
 }
  \ee
and the subset $\{\psi_k^*\}$ is complete
 in
$L^2_{\rho^*}(\re)$.

\end{lemma}

As usual, if $\{\psi_k\}$ is the adjoint (Riesz) basis of
eigenfunctions of the adjoint operator
 \be
 \label{BBB}
 \tex{
 \BB= - D_y^4 + \frac 14\, y D_y + \frac 14\, I \inB
 L^2_\rho(\re) \withA \rho= \frac 1{\rho^*},
 }
\ee
   with the same spectrum \ef{SpecN}, the bi-orthonormality
condition holds in $L^2(\re)$:
 \be
 \label{Ortog}
\langle \psi_k, \psi_l^* \rangle = \d_{k l} \quad \mbox{for any}
\quad k, \,\, l.
 \ee
 The generating formula for $\psi_k$ is as follows:
  \be
  \label{gen11}
   \tex{
    \psi_k(y)= \frac{(-1)^k}{ \sqrt{k !}}\, F^{(k)}(y), \quad k
    =0,1,2,...\, ,
    }
    \ee
    where $F(y)$ is the rescaled kernel of the fundamental
    solution $b(x,t)$ of  $D_t+D_x^4$, i.e.,
    \be
    \label{gen12}
     \tex{
     b(x,t)=t^{-\frac 14} F(y), \,\, y= \frac x{t^{1/4}}; \,\,\,
     b_t=-b_{xxxx}, \,\, b(x,0)=\d(x) \, \Longrightarrow \, \BB F=0, \,\, \int F=1.
 }
 \ee

\subsection{Why self-similar blow-up is expected to be generic}

This can be connected with the spectral gaps in
$\s(\BB^*+I)=\{1-\frac k4\}$ of the linearized operator in
\ef{4.1}. Indeed, we observe here four positive unstable modes:
 \be
 \label{ss1}
 k=0 \andA 1, \andA k=2 \andA 3.
  \ee
  The first group is not taken into account since correspond to
  standard instabilities relative to perturbations of the time $T=0$
  and the space $x=0$ of the blow-up point.

The crucial is the actual unstable mode with
 \be
 \label{un1}
  \tex{
 k=2 \whereA \l_2=\frac 12 \andA \psi_2^*(y)= \frac 1{\sqrt 2} \,
 y^2.
  }
  \ee
 We claim that this unstable mode ensures  convergence to a
 non-trivial similarity stationary profile, which we used to
 denote by $f_0(y)$; see extra explanations in \cite[\S~3]{GW1}. There is no still a proof of this
 intriguing fact on such a heteroclinic connection $f_* \to f_0$, but it sounds rather reliably.
 Hence, the stable manifold of $f_0$ includes a part of the
 unstable manifold of the constant equilibrium $f_*=(p-1)^{-1/(p-1)}$,
 so that linearized patterns to be constructed must be less
 stable than $f_0$, which then is expected to represent   a most structurally
 stable (generic) blow-up pattern.

Note that the third unstable mode in \ef{ss1} with $k=3$ has
nothing to do with the second similarity profile $f_1$. Indeed,
this mode corresponds to odd (anti-symmetric) perturbations, while
$f_1(y)$ is an even function. The origin of existence of $f_1$
 is more subtle and is explained in \cite{BGW1} by branching and homotopy theory.
Moreover, according to \cite[\S~3.2]{GW1}, $f_1(y)$ is connected
with the trivial equilibrium $0$ via a semi-stable centre manifold
behaviour to be presented below.

\subsection{Inner expansion}

Thus, in the {\em Inner Region} characterized by compact subsets
in the similarity variable $y$, we assume a centre or a stable
subspace behaviour as $\t \to +\iy$ for the linearized operator
$\BB^*+I$ (see \cite{GalCr} for rigorous details in a simpler
related problem):
 \be
 \label{4.5}
  \begin{matrix}
 \mbox{centre:} \quad Y(y,\t)=a(\t) \psi_4^*(y) + w^\bot
 \quad (\l_4=-1, \,\, k=4),\quad\,\,\ssk\ssk\\
 \mbox{stable:} \quad  Y(y,\t) =-C {\mathrm e}^{\l_k \t} \psi_k^*(y)+w^\bot
 \quad (\l_k <-1, \,\, k>4).
    \end{matrix}
    \ee
 As usual for the classic equation \ef{m1}, we restrict to even
 values only,
   \be
   \label{par1}
  k=4,6,8,...\, ,
   \ee
  since existence of patterns with odd $k$'s, i.e., having a
  non-symmetric blow-up structures is rather suspicious. However,
  for higher-order diffusion, such patterns cannot be excluded
  entirely, though this is not our business here.

For the centre subspace behaviour in \ef{4.5}, substituting the
eigenfunctions expansion into equation \ef{4.1} yields the
following coefficient:
 \be
 \label{4.6}
  \tex{
  \dot a=a^2 \g_0+... \whereA \g_0= c_0 \langle (\psi_4^*)^2,
\psi_4 \rangle
  \LongA a(\t) =- \frac {1}{\g_0\, \t}+... \, .
  }
  \ee
 It is crucial that \ef{4.6} shows a clear {\em semi-stable}
(saddle-node) structure of the equilibrium $a=0$, which
essentially depends on the sign of the coefficient $\g_0$; see
below.

 Note that for the matching purposes, one needs to
 assume that (see details in \cite{Gal2m, GW1}):
 \be
 \label{4.7}
\mbox{if $\psi_k^*(0)>0$, then} \quad  \g_0>0 \andA C>0.
  \ee
Actually, for $k=4$ it is calculated explicitly that
 \be
 \label{4.71d}
  \tex{
  \g_0=-c_0 136 \sqrt 6<0,
   }
   \ee
 so that the centre manifold patterns with the positive eigenfunction
 \be
 \label{4.72d}
  \tex{
   \psi_4^*(y)= \frac 1{\sqrt{24}}\, (y^4+ 24)
    }
    \ee
correspond to solutions that blow-up on finite interfaces,
\cite[\S~3]{GW1}. It seems reasonable that a full justification of
such a behaviour can be done along the lines of classic invariant
manifold theory (see e.g. \cite{Lun}), though can be very
difficult.

   Overall, the whole variety of such asymptotics is
characterized as follows:
 \be
 \label{YY1}
  Y(y,\t) = a(\t) y^k +...\, \whereA a(\t) = -
 \left\{
  \begin{matrix}
 \,\,\,
  \frac
  1{\g_0 \t}+... \,\,\,
  \quad
  \mbox{for}\quad k=4,\ssk\\
  C{\mathrm e}^{\l_k \t}+... \quad \mbox{for} \quad k>4.
 \end{matrix}
 \right.
  \ee

\subsection{Outer region: matching}

 The asymptotics
\ef{4.5} is known  \cite{Gal2m} to admit matching with the {\em
Outer Region}, being a Hamilton--Jacobi (H--J) one. More
precisely, in the centre case with $k=4$, according to \ef{4.5},
\ef{4.6}, we introduce the outer variable and obtain from \ef{3.2}
the following perturbed H--J equation:
 \be
 \label{4.71}
 \tex{
 \xi= \frac y {\t^{1/4}} \LongA v_\t= - \frac 1 4\, \xi v_\xi - \frac 1{p-1}\, v + |v|^{p-1}
 v + \frac 1 \t \, \big( \frac 14\, \xi  v_\xi - v_{\xi\xi\xi\xi} \big).
 }
 \ee
Passing to the limit as $\t \to +\iy$ in such singularly perturbed
PDEs is not easy at all even in the second-order case (see a
number of various  applications in \cite{AMGV}). Though currently
not rigorously, 
 we assume
stabilization to the stationary solutions $f(\xi)$ satisfying the
unperturbed H--J equation:
 \be
 \label{4.72}
 \tex{
 - \frac 14\, \xi f'- \frac 1{p-1}\, f + |f|^{p-1}
 f=0 \inB \re.
 }
 \ee
 This is solved via characteristics, where we
 have to choose the solution satisfying
 \ef{YY1}:
  \be
  \label{4.73}
   \tex{
  f(\xi)= (p-1)^{-\frac 1{p-1}}  - \frac 1{\g_0}\, \xi^4
  +... \asA
  \xi \to 0.
  }
  \ee
 Hence, integrating \ef{4.72} with the condition \ef{4.73} yields
 the unique H--J profile
  \be
  \label{f1}
  \tex{
  f_0(\xi)=f_* \, (1+c_* \xi^4)^{-\frac 1{p-1}} \whereA c_* =\frac
  1{\g_0}\, (p-1)^{\frac p{p-1}}.
   }
   \ee
   Note that since $\g_0<0$ according to \ef{4.71d}, the profile
   \ef{f1} blows up as $c_* \xi^4 \to -1$; see analogous phenomena
   in \cite{GW1}.
    For the $2m$th-order PDE \ef{mm98} with odd
   $m=3,5,...$\, , we have $\g_0>0$ and then \ef{f1} can represent
   standard blow-up patterns, \cite{Gal2m}.

Analogously, for the stable behaviour for $k>4$ in \ef{4.5}, we
use the following change:
\be
 \label{4.71S}
 \tex{
 \xi=  {\mathrm e}^{ \frac{4-k}{4 k}\, \t}\, y  \LongA v_\t=
 - \frac 1{k}\, \xi  v_\xi- \frac 1{p-1}\, v + |v|^{p-1}
 v - {\mathrm e}^{- \frac {k-4}{k}\, \t} \,  v_{\xi\xi\xi\xi}.
 }
 \ee
Passage to the limit $\t \to +\iy$ and matching with Inner Region
are analogous. The H--J profile is as follows:
 \be
  \label{f2}
  \tex{
  f_C(\xi)=f_* \, (1+c\, \xi^k)^{-\frac 1{p-1}} \whereA c=
  C\, (p-1)^{\frac p{p-1}}>0.
   }
   \ee

Overall,  according to matching conditions \ef{YY1} and
\ef{4.71S},
the whole set of  such blow-up patterns
consists of  a countable
 set of different solutions for $k=6,8,...$ (as we have mentioned,
 the case $k=4$ works for $m$ odd only, and then $k=2m$).



\subsection{Final time profiles: the rout to complete blow-up}

As a natural counterpart of the above asymptotic analysis, we
introduce the {\em final time profiles} of solutions
\cite[\S~4]{Gal2m}:

\ssk

\noi{\sc Centre subspace patterns: $k=4$.} This  analysis is
formal, since by \ef{4.71d} and \ef{f1} with $c_*<0$  such
patterns are not bounded. However, we keep the derivation formally
assuming that  case $\g_0>0$, which applies to equations \ef{mm98}
with odd $m=3,5,7,...$\,, \cite{Gal2m}. Thus,
 marching of the outer region with the Outer Region II for $|x|>0$
 small, yields the following behaviour:
  \be
  \label{ft1}
   \tex{
u(x,0^-) = C_* |x|^{-\frac{4}{p-1}}\, \big|\ln|x|\,\big|^{\frac
1{p-1}}(1+o(1)),
 }
 \ee
where $C_*>0$ is a constant depending on the parameters $p$ and
$\g_0$ (assumed here to be positive) only and hence is independent
of initial data $u_0$. A more careful passing to the limit as $t
\to 0^-$ in \ef{3.1MM} by using the limit profile  \ef{f1} shows
that
 \be
 \label{f33}
  \tex{
  C_*= f_* \, \big(\frac 4{c_*}\big)^{\frac 1{p-1}}.
   }
   \ee


\noi{\sc Stable subspace patterns: $k>4$.} Such stable subspace
patterns can be truly constructed for all $m=2,3,...$\,. Matching
with the Outer Region II yields, as $x \to 0$,
 \be
 \label{ft2}
 u(x,0^-) =  C_*
|x|^{-\frac{k}{p-1}}(1+o(1)) \forA k=6,8,...\, ,
 \ee
 where constants  $C_*$ depend on initial function $u_0$ through the earlier constant
 $C$ in \ef{YY1}.

\ssk

Comparing the self-similar profile \ef{dd3} with those in \ef{ft1}
and \ef{ft2} shows the actual origin of {\em complete blow-up}:
the later ones contain stronger singularities as initial data at
$x=0$ posed at $t=0$. This is easier to see for the case $k=4$
(actually, nonexistent by \ef{4.71d}), where the profiles differ
by the unbounded factor
 \be
 \label{fac1}
  \tex{
  \sim \big|\ln|x|\,\big|^{\frac
1{p-1}} \to +\iy \asA x \to 0.
 }
 \ee
Surely, the same and in a stronger manner happens for
$k=6,8,...\,$.

 Bearing in mind our self-similar blow-up patterns \ef{2.1} with
 the asymptotics \ef{dd1}, the unbounded factor in \ef{fac1}
 actually means that, in a certain natural but formal sense, we
 have to look for a similarity pattern with
  \be
  \label{fac2}
   C_1=+\iy.
   \ee
 Moreover, for existence of a regular extension for $t>0$, this
would demand existence of a global patters \ef{3.1}, \ef{3dd1}
also with the coefficient \ef{fac2}. Of course, this is not
possible that somehow reflects the nonexistence of a proper
extension, meaning complete blow-up.

The actual proof of nonexistence of a regular solution of \ef{m2}
with data \ef{ft2} is not easy. This assumes proving nonexistence
of a corresponding {\em very singular solution} for \ef{m2}, which
is a well-known issue for second-order semilinear and quasilinear
heat equations. However, such results are also known for
higher-order nonlinear parabolic equations with various
nonlinearities, where the analysis without the Maximum Principle
becomes essentially more involved. We refer to functional methods
in \cite{GSVSS1}, where necessary earlier references on the
subject can be found.

 \section{Final remark: on evolution completeness of blow-up
 patterns}
  \label{S.5}

 After introducing and discussing all the blow-up patterns, we
 finally are in a position to announce another important aspect
 of our study. Namely, we claim that the whole set of self-similar
 and linearized patterns introduced above are {\em evolutionary
 complete} in a natural sense; see \cite{CompG}. This purely  means that
 the self-similar blow-up patterns \ef{2.3} together with all the
 linearized ones in Section \ref{S4} are expected to describe all
 possible types of blow-up that can occur for the PDE \ef{m2}. We
 connect this claim with the standard completeness of the
 eigenfunctions \ef{psidec} in $L^2_{\rho^*}$. However, since by
 \ef{ss1} we have thrown away two unstable patterns for $k=2$ and
 $k=3$ (the modes with $k=0$ and 1 do not occur in the evolution),
 one needs to replace these two ones. And this is done by
 including precisely {\em two} nonlinear self-similar patterns
 \ef{2.1} with profiles $f_0$ and $f_1$.
 Note that this may look irrelevant since the mode with $k=3$ is
 spatially
 odd, while $f_1(y)$ was shown to be even. In fact, this is not
 that important, since we always mean a certain symmetrization of
 blow-up behaviour as $t \to T^-$ (no proof is still available),
  so that including two ``nonlinear
 eigenfunctions" $f_{0,1}(y)$ is expected to be enough
 to govern non-trivial self-similar blow-up.

What is more crucial is that a standard centre manifold behaviour
is prohibited by the sign in \ef{4.71d}. Therefore, finally,
according to evolution completeness arguments:
 $$
  \fbox{$
 \mbox{two similarity profiles $f_{0,1}(y)$ replace two modes,
 with $k=2$ and 4.}
 $}
 $$
Here we observe a usual ``dimension preservation" by including
into the countable set of linearized patterns (minus two modes for
$k=2,4$) two nonlinear eigenfunctions $f_{0,1}(y)$.

Of course, another and obviously
 inevitable key aspect of the
completeness speculation above is {\em nonexistence} of Type II
blow-up in \ef{S5} for \ef{m2}. Indeed, Type II blow-up patterns
would destroy any completeness of the above functional set. Note
that such Type II patterns were obtained for \ef{r3} for some
large $p$ in the supercritical Sobolev range
 $$
 \tex{
 p \ge  p_{\rm Sob}= \frac{N+4}{N-4} \quad (N>4);
 }
 $$
see \cite[\S~5,\,6]{GalBlow5} for details. We have some reliable
evidence that such patterns are nonexistent in the subcritical
range $p< p_{\rm Sob}$ and, in particular, for $N=1$, i.e., for
\ef{m2}. The point is that, in the natural rescaled sense, such
patterns on smaller compact subsets around the point $\{0,0^-\}$
must be governed by a  regular stationary solution $W$ satisfying
  \be
  \label{W11}
 -\D^2 W + |W|^{p-1}W=0, \quad W(0)=1.
  \ee
In the subcritical range, it is guaranteed that radial solutions
of \ef{W11} are highly oscillatory \cite{Gaz06}, which makes very
unlikely to match such a behaviour with the outer region to get an
acceptable blow-up pattern of \ef{m2} in the
$\{x,t,u\}$-variables. We again refer to \cite[\S~5]{GalBlow5},
and must admit that a rigorous proof of such a nonexistence is
absent.





\begin{thebibliography} {10}


 \bibitem
 {AL84}
M.M.~Ad'jutov and L.A.~Lepin, {Absence of blowing up similarity
structures in a medium with a source for constant thermal
conductivity}, Differ. Equat.,
 {\bf 20} (1984), 1279--1281.

 \bibitem
  {Barb}
G.~Barbatis, {\em Explicit estimates on the fundamental solution
  of higher-order parabolic equations with measurable coefficients}, J.~Differ.
   Equat., {\bf 174} (2001), 442--463.


    \bibitem
  {Barb04}
G.~Barbatis, {\em Sharp heat-kernel estimates for
 higher-order operators with singular coefficients},
 Proc. Edinb. Math. Soc. (2),
  {\bf 47} (2004), 53--67.




\bibitem 
{BebEb} J.~Bebernes and D.~Eberly, {Mathematical Problems in
Combustion Theory}, Appl.
  Math. Sci., Vol. {\bf 83}, Springer-Verlag, Berlin, 1989.




 \bibitem{BGW1}
C.J.~Budd, V.A.~Galaktionov, and J.F.~Williams, \emph{Self-similar
blow-up
  in higher-order semilinear parabolic equations}, SIAM J.~Appl. Math.,
 {\bf 64} (2004), 1775--1809.


\bibitem{CP}
C.J.~Chapman and M.R.E.~Proctor, \emph{Nonlinear
{R}ayleigh-{B}enard convection
  between poorly conducting boundaries}, J.~Fluid Mech., \textbf{101} (1980),
  759--782.


\bibitem
{CodL} E.A.~Coddington and N.~Levinson, {Theory of Ordinary
Differential Equations}, McGraw-Hill Book Company, Inc., New
York/London, 1955.







\bibitem{Eg4}
Yu.V.~Egorov, V.A.~Galaktionov, V.A.~Kondratiev, and
S.I.~Pohozaev,
  {\em Global solutions of higher-order semilinear
  parabolic equations in the supercritical range}, {Adv. Differ.
  Equat.,}
{\bf 9} (2004), 1009--1038.


     \bibitem{Bl4}
J.D.~Evans, V.A.~Galaktionov, and J.R.~King, \emph{Blow-up
similarity solutions of  the fourth-order unstable thin film
equation},  Euro J.~Appl. Math., {\bf 18} (2007), 195--231.


\bibitem
  {EGW1}
J.D.~Evans, V.A.~Galaktionov, and J.F.~Williams, {\em Blow-up and
global asymptotics of  the limit unstable Cahn-Hilliard equation},
SIAM J.~Math. Anal., {\bf 38} (2006), 64--102.


\bibitem
{Fr-K}
D.A.~Frank-Kamenetskii, \emph{Towards temperature distributions in
a reaction
  vessel and the stationary theory of thermal explosion}, Doklady Acad. Nauk
  SSSR \textbf{18} (1938), 411--412.

\bibitem{F-K}
D.A.~Frank-Kamenetskii, {Diffusion and {H}eat {T}ransfer in
{C}hemical {K}inetics},
  Plenum Press, New York, 1969.


\bibitem 
{Gal2m} V.A. Galaktionov, {\em On a spectrum of blow-up patterns
for a higher-order semilinear parabolic equations}, Proc. Royal
Soc. London A, {\bf 457} (2001), 1-21.




 \bibitem
 {GalCr}
  V.A.~Galaktionov, {\em Critical global asymptotics in
  higher-order semilinear parabolic equations},
 Int.~J. Math. Math. Sci., {\bf 60} (2003),  3809--3825.


\bibitem
 {CompG}
 V.A.~Galaktionov,
{\em Evolution completeness of separable solutions of non-linear
diffusion equations in bounded domains}, Math. Meth. Appl. Sci.,
{\bf 27} (2004), 1755--1770.




\bibitem 
 {GalGeom}
 V.A.~Galaktionov, {\rm Geometric Sturmian  Theory of Nonlinear
 Parabolic Equations and Applications}, Chapman$\,\&\,$Hall/CRC, Boca Raton,
Florida,
 2004.


\bibitem
 {GalJMP}
  V.A.~Galaktionov,
 {\em On blow-up  space jets for the Navier--Stokes equations in $\re^3$
 with convergence to Euler equations},  J.~Math. Phys., {\bf 49}
 (2008), 113101.




\bibitem
 {GalBlow5}
 V.A.~Galaktionov,
 {\em Five types of blow-up in a semilinear fourth-order reaction-diffusion equation:
 an analytic-numerical approach}, Nonlinearity,
   submitted (arXiv:0901.4307).



\bibitem
 {GHUni}
 V.A.~Galaktinov and P.J.~Harwin,
  {\em Non-uniqueness and  global similarity solutions   for
a higher-order semilinear parabolic equation}, Nonlinearity, {\bf
18} (2005), 717--746.


\bibitem{GP5}
V.A.~Galaktionov and S.A.~Posashkov, {\em Application of
  a new comparison theorem for unbounded solutions of nonlinear
  parabolic equations}, Differ. Equat.,
  {\bf 22} (1986), 809--815.


\bibitem 
{GSVSS1}
 V.A.~Galaktionov and A.E.~Shishkov, {\em  Higher-order
quasilinear parabolic equations
 with singular initial data}, Comm. Contemp. Math., {\bf 8} (2006), 331--354.



\bibitem
 {GSVR} V.A.~Galaktionov and S.R.~Svirshchevskii, Exact Solutions and
 Invariant Subspaces of Nonlinear Partial Differential Equations in Mechanics and Physics,
  Chapman$\,\&\,$Hall/CRC, Boca Raton,
Florida,
 2007.



  \bibitem 
{GV}   V.A.~Galaktionov and J.L.~Vazquez, {\em Continuation of
blow-up solutions of nonlinear heat equations in several space
dimensions,} { Comm. Pure Appl. Math.,} {\bf 50} (1997), 1-68.



\bibitem 
{AMGV}
 V.A.~Galaktionov and J.L.~Vazquez, {A Stability Technique  for Evolution
Partial Differential Equations.
 A Dynamical Systems Approach},
{Progr. in Nonl. Differ.
  Equat. and Their Appl.,} Vol. {\bf 56},
Birkh\"auser, Boston/Berlin, 2004.


\bibitem
{GW1}
V.A.~Galaktionov and J.F.~Williams, \emph{Blow-up in a
fourth-order
  semilinear parabolc equation from explosion-convection theory},
  Euro~J.~Appl. Math., {\bf 14} (2003), 745--764.


  \bibitem
{Gaz06}
 F.~Gazzola and H.-C.~Grunau,  {\em Radial entire solutions for
 supercritical biharmonic equations},
 Math. Ann.,
 {\bf 334} (2006), 905--936.


\bibitem{GS}
V.L.~Gertsberg and G.I.~Sivashinsky, \emph{Large cells in
nonlinear
  rayleigh-benard convection}, Prog. Theor. Phys., \textbf{66} (1981),
  1219--1229.


 \bibitem
  {GK85}
 Y.~Giga and R.V.~Kohn, {\em Asymptotically self-similar blow-up of
 semilinear heat equations}, Comm. Pure Appl. Math., {\bf 38}
 (1985), 297--319.

 \bibitem
   {Ham95}
  R.~Hamilton, {\em The formation of singularities in the Riccu
  flow. Surveys in Differ. Geom.,} Vol. II (Cambridge, MA, 1993),
  pp. 7-136, Int. Press, Cambridge, MA, 1995.


\bibitem{JMS}
G.~Joulin, A.B. Mikishev, and G.I. Sivashinsky, \emph{A
  {S}emenov-{R}ayleigh-{B}enard problem}, Preprint.



\bibitem
 {Ler34}
  J.~Leray, {\em Sur le mouvement d'un liquide visqueux emplissant l'espace},
   Acta Math., {\bf 63}
  (1934), 193--248.



\bibitem
{Lun} A.~Lunardi, {\rm Analytic Semigroups and Optimal Regularity
in Parabolic Problems}, Birkh\"auser, Basel/Berlin, 1995.


\bibitem 
 {MitPoh}
E.~Mitidieri and S.I.~Pohozaev, {\rm Apriori Estimates and Blow-up
of Solutions to Nonlinear Partial Differential Equations and
Inequalities},  Proc. Steklov Inst. Math., Vol. {\bf 234}, Intern.
Acad. Publ. Comp. Nauka/Interperiodica, Moscow, 2001.

%






\bibitem
 {Pao}
 C.V.~Pao,  {\rm Nonlinear Parabolic and Elliptic Equations},
  Plenum Press, New York, 1992.  



\bibitem
 {Plan07}
F.~Planchon  and P.~Rapha\"el, {\em Existence and stability of the
log--log blow-up dynamics for the $L^2$-critical  nonlinear
Schr\"odinger equation in a domain}, Ann. Henri Poincar\'e, {\bf
8} (2007), 1177--1219.


\bibitem
{PelTroy}
  L.A.~Peletier and W.C.~Troy, {\rm Spatial Patterns.
Higher Order Models in Physics and Mechanics}, Birkh\"auser,
Boston/Berlin, 2001.




\bibitem
 {QSupl}
  P.~Quittner and P.~Souplet, {\rm Superlinear Parabolic Problems
  and Their Equilibria}, Birkh\"auser, 2007.





\bibitem 
{SGKM}  A.A. Samarskii, V.A. Galaktionov, S.P. Kurdyumov, and A.P.
Mikhailov, {Blow-up in Quasilinear Parabolic Equations}, \rm
Walter de Gruyter, Berlin/New York, 1995.

\bibitem{Sem}
N.~Semenov, {Chemical {K}inetics and {C}hain {R}eaction},
Clarendon Press,
  Oxford, 1935.



\bibitem 
{VainbergTr} M.A.~Vainberg and V.A.~Trenogin, {\rm Theory of
Branching of Solutions of Non-Linear Equations}, Noordhoff Int.
Publ., Leiden, 1974.



  \bibitem 
{ZBLM}
 Ya.B.~Zel'dovich, G.I.~Barenblatt, V.B.~Librovich,
and G.M.~Makhviladze, {\rm The Mathematical Theory of Combustion
and Explosions}, Consultants Bureau [Plenum], New York, 1985.




\end{thebibliography}
 \end{document}